\documentclass[12pt]{amsart}
\usepackage{amssymb}
\usepackage{amsmath}
\usepackage{amsthm}

\DeclareMathOperator{\ad}{ad}
\DeclareMathOperator{\id}{id}

\renewcommand{\c}{{\mathcal C}}

\DeclareMathOperator{\Hom}{Hom}

\newcommand{\g}{{\mathfrak g}}
\newcommand{\I}{{\mathcal I}}
\newcommand{\J}{{\mathcal J}}

\renewcommand{\ll}{{\langle\langle}}
\renewcommand{\L}{\Lambda}

\newcommand{\rr}{{\rangle\rangle}}
\newcommand{\N}{{\mathcal N}}

\newcommand{\C}{\ensuremath{\mathbb{C}}}
\newcommand{\R}{\ensuremath{\mathbb{R}}}

\newcommand{\p}{\partial}
\newcommand{\s}{{\rm Symb}}
\renewcommand{\S}{\Sigma}
\newcommand{\T}{{T^\ast M}}
\newcommand{\D}{{\mathcal D}}
\newcommand{\E}{{\mathcal E}}
\newcommand{\U}{{\mathcal U}}

\newtheorem{lemma}{Lemma}
\newtheorem{proposition}{Proposition}
\newtheorem{theorem}{Theorem}

\begin{document}
\title[Formal symplectic groupoid]
{Formal symplectic groupoid of a deformation quantization}
\author[A.V. Karabegov]{Alexander V. Karabegov}\footnote{Research was partially supported by an ACU Math/Science grant.}
\address[Alexander V. Karabegov]{Department of Mathematics and Computer Science, Abilene Christian University, ACU Box 28012, 252 Foster Science Building, Abilene, TX 79699-8012}
\email{alexander.karabegov@math.acu.edu}

\begin{abstract}
We give a self-contained algebraic description of a formal symplectic groupoid over a Poisson manifold $M$. To each natural star product on $M$ we then associate a canonical formal symplectic groupoid over $M$. Finally, we construct a unique formal symplectic groupoid `with separation of variables' over an arbitrary K\"ahler-Poisson manifold.
\end{abstract}

\date{June 19, 2004}
\maketitle

\section{Introduction}

Symplectic groupoids are semiclassical geometric objects whose heu\-ristic quantum counterparts are associative algebras treated as quantum objects. In \cite{CDF}, \cite{Deq}, and \cite{Inv} evidence was given that the star algebra of a deformation quantization gives rise to a formal analogue of a symplectic groupoid. In this paper we give a global definition of a formal symplectic groupoid and show that to each natural deformation quantization (in the sense of Gutt and Rawnsley, \cite{GR}) there corresponds a canonical formal symplectic groupoid.

Symplectic groupoids were introduced independently by Karas\"ev \cite{Ka}, Weinstein \cite{W}, and Zakrzewski \cite{Z}. Recall that a local symplectic groupoid is an object that has the properties of a symplectic groupoid in which the multiplication is local, being only defined in a neighborhood of the unit space. It was proved in \cite{Ka} and \cite{W} that for any Poisson manifold $M$ there exists a local symplectic groupoid over $M$ that `integrates' it.  In \cite{CDF} A. S. Cattaneo, B. Dherin, and G. Felder considered the formal integration problem for $\R^n$ endowed with an arbitrary Poisson structure, whose solution is given by a formal symplectic groupoid. They start with the zero Poisson structure on $\R^n$, the corresponding trivial symplectic groupoid $T^* \R^n$, and a generating function of the Lagrangian product space of this groupoid. A formal symplectic groupoid is then defined in terms of a formal deformation of that trivial generating function. One of the main results of \cite{CDF} is an explicit formula for a generating function that delivers the formal symplectic groupoid related to the Kontsevich star product. The approach to formal symplectic groupoids developed in \cite{CDF} demonstrates the relationship between geometric and algebraic deformations described in \cite{WTG}.

One can take an alternative approach to the definition of a formal symplectic groupoid 
over a Poisson manifold $M$ (which leads to the same object) by replacing the symplectic manifold $\S$ on which a (local) symplectic groupoid over $M$ is defined,  with the formal neighborhood $(\S,\L)$ of its unit space $\L$ (see the definition of a formal neighborhood in Section \ref{S:def}). We use a simple model of the algebra of formal functions on $(\S,\L)$ which is reminiscent of the Hopf algebroid constructed by Vainerman in \cite{V}. This model provides effective means to check the axioms of a formal symplectic groupoid and to do the calculations.

In Section \ref{S:def} we state formal analogues of the axioms of a symplectic groupoid and give a definition of a formal symplectic groupoid over a given Poisson manifold $M$. Such a formal groupoid is defined on the formal neighborhood of a Lagrangian submanifold of a symplectic manifold. In Section \ref{S:formsymp} we give a self-contained algebraic description of a formal symplectic groupoid and show that a strict formal symplectic realization of an arbitrary Poisson manifold $M$ gives rise to a unique formal symplectic groupoid over $M$ whose source mapping is given by that formal symplectic realization. In Section \ref{S:isom} we describe the space of all formal symplectic groupoids over $M$ which are defined on a given formal neighborhood of a Lagrangian submanifold of a symplectic manifold. In Section \ref{S:cannat} we relate to each natural deformation quantization on $M$ a canonical formal symplectic groupoid. In Section \ref{S:defquant} we prove that any deformation quantization with separation of variables on a K\"ahler-Poisson manifold $M$ is natural and show that its canonical formal symplectic groupoid has a property which we call `separation of variables'. Finally, in Section \ref{S:formgroupsv} we prove that for an arbitrary K\"ahler-Poisson manifold $M$ there exists a unique formal symplectic groupoid with separation of variables over $M$.

\section{Definition of a formal symplectic groupoid}\label{S:def}

A symplectic groupoid over a Poisson manifold $(M, \{\cdot,\cdot\}_M)$ is a symplectic manifold $\S$ endowed with the associated Poisson source mapping $s: \S \to M$, the anti-Poisson target mapping $t: \S \to M$, which both are surjective submersions, the antisymplectic involutive inverse mapping $i: \S \to \S$, and the unit mapping $\epsilon: M \to \S$, which is an embedding. The image $\L = \epsilon(M)$ of the unit mapping is the Lagrangian unit space of the symplectic groupoid. Denote by $\S^n$ the Cartesian product of $n$ copies of the manifold $\S$ and by $\S_n$ the submanifold of $\S^n$ formed by the $n$-tuples $(\alpha_1, \ldots, \alpha_n) \in \S^n$ such that $t(\alpha_k) = s(\alpha_{k+1}), \ 1 \leq k \leq n-1$. The coisotropic submanifold $\S_2$ of $\S \times \S $ is the domain of the groupoid multiplication $m: \S_2 \to \S$. For $\alpha,\beta \in \S_2$ we write $m(\alpha, \beta) = \alpha\beta$. 
The groupoid multiplication is associative. For $(\alpha,\beta,\gamma) \in \S_3$ the associativity condition $\alpha(\beta\gamma) = (\alpha\beta)\gamma$ holds. The graph $\Gamma = \{(\alpha, \beta, \gamma) \, | \, (\beta, \gamma) \in \S_2, \alpha = \beta\gamma\}$ of the groupoid multiplication (the product space) is a Lagrangian submanifold of $\S \times \bar \S \times \bar \S$, where $\bar \S$ is a copy of the manifold $\S$ endowed with the opposite symplectic structure. 

The groupoid operations satisfy the following axioms. For any composable $\alpha, \beta \in \S$ and $x \in M$
\begin{align*}
(A1) \ s(\alpha\beta) = s(\alpha), \ (A2) \ t(\alpha\beta) = t(\beta), \ (A3) \ s \circ \epsilon = \id_M, \\
 (A4) \ t \circ \epsilon = \id_M, \ (A5) \ \epsilon (s(\alpha))\alpha = \alpha, (A6) \ \alpha \epsilon (t (\alpha)) = \alpha, \\ 
(A7) \ s(i(\alpha)) = t(\alpha), (A8) \ \alpha i(\alpha) = \epsilon (s (\alpha)),  \ (A9) \  i (\alpha)\alpha = \epsilon (t (\alpha)).
\end{align*}

Recall the definition of the formal neighborhood $(X,Y)$ of a submanifold $Y$ of a manifold $X$. Let $Y$ be a closed $k$-dimensional submanifold of a real $n$-dimensional manifold $X$ and $I_Y \subset C^\infty(X)$ be the ideal of smooth functions on $X$ vanishing on $Y$. Then the quotient algebra $C^\infty(X,Y) := C^\infty(X)/I_Y^\infty$, where $I_Y^\infty = \cap_{l = 1}^\infty I_Y^l$, can be thought of as the algebra of smooth functions on the formal neighborhood $(X,Y)$ of the submanifold $Y$ in $X$. If $U \subset X$ is a local coordinate chart on $X$ with coordinates $\{x^i\}$ such that $U\cap Y$ is given by the equations $x^{k+1} = 0,\ldots,x^n =  0$, then $C^\infty(U,U \cap Y)$ is isomorphic to $C^\infty(U\cap Y)[[x^{k+1},\ldots, x^n]]$, where the isomorphism is established via the formal Taylor expansion of the functions on $U$ in the variables $x^{k+1},\ldots, x^n$. Thus the formal neighborhood $(X,Y)$ of $Y \subset X$ is the ringed space $Y$ with the sheaf of rings whose global sections form the algebra $C^\infty(X,Y)$. Let $Y_i$ be a submanifold of a manifold $X_i$ for $i = 1,2$. If $f: X_1 \to X_2$ is a mapping such that $f(Y_1) \subset Y_2$, then $f^*(I_{Y_2}) \subset I_{Y_1}$. Therefore the mapping $f$ induces the dual morphism of algebras $f^*: C^\infty(X_2,Y_2) \to C^\infty(X_1,Y_1)$.

Denote by $\L^n$ the Cartesian product of $n$ copies of the manifold $\L$ and by $\L_n$ the diagonal of $\L^n$. Notice that $\S_n \cap \L^n = \L_n$. The algebra $C^\infty(\S)$ is a Poisson algebra with respect to the natural Poisson bracket $\{\cdot,\cdot\}_\S$ on $\S$. The space $C^\infty(\S,\L)$ inherits a structure of  Poisson algebra  from $C^\infty(\S)$. We will use the same notation $\{\cdot,\cdot\}_\S$ for the induced Poisson bracket on $C^\infty(\S,\L)$. Similarly, denote by $\{\cdot,\cdot\}_{\S^n}$ the Poisson bracket on $C^\infty(\S^n)$ corresponding to the product Poisson structure, and the induced bracket on $C^\infty(\S^n, \L^n)$. Let $\iota : \S_2 \to \S \times \S$ be the inclusion mapping. We will say that functions $F\in C^\infty(\S)$ and $G\in C^\infty(\S\times \S)$ such that $m^*F = \iota^* G$ agree on $\S_2$. The functions $F\in C^\infty(\S)$ and $G\in C^\infty(\S\times \S)$ agree on $\S_2$ if and only if the function $F \otimes 1 \otimes 1 - 1 \otimes G \in C^\infty(\S \times \bar \S \times \bar \S)$ vanishes on the product space $\Gamma$. Since $\Gamma$ is a Lagrangian submanifold of 
$\S \times \bar \S \times \bar \S$, the Poisson bracket of two functions vanishing on $\Gamma$ also vanishes on $\Gamma$. For functions $F_i\in C^\infty(\S)$ and $G_i\in C^\infty(\S\times \S),\ i = 1,2,$ the Poisson bracket of $F_1 \otimes 1 \otimes 1 - 1 \otimes G_1$ and $F_2 \otimes 1 \otimes 1 - 1 \otimes G_2$ equals
\[
   \{F_1,F_2\}_\S \otimes 1\otimes 1 - 1 \otimes \{G_1,G_2\}_{\S^2},
\]
whence we obtain the following lemma.
\begin{lemma}\label{L:propp}
If functions $F_i\in C^\infty(\S)$ and $G_i\in C^\infty(\S\times \S),\ i = 1,2,$ agree on $\S_2$, then the Poisson brackets $\{F_1,F_2\}_\S$ and $\{G_1,G_2\}_{\S^2}$ also agree on $\S_2$.
\end{lemma}

The multiplication $m : \S_2 \to \S$ identifies $\L_2$ with $\L$ and thus induces the comultiplication mapping
\[
m^* : C^\infty(\S,\L) \to C^\infty(\S_2,\L_2).
\] 
Denote by $\iota_n : \S_n \to \S^n$ the inclusion mapping. In particular, $\iota = \iota_2$. 
Since the mapping $\iota_n$ maps $\L_n$ to $\L^n$, it induces the algebra morphism
\[
         \iota^*_n : C^\infty(\S^n,\L^n) \to C^\infty(\S_n,\L_n).
\]
We say that elements $F\in C^\infty(\S,\L)$ and  $G\in C^\infty(\S^2,\L^2)$ agree on $C^\infty(\S_2,\L_2)$ if $m^*F = \iota^* G$ in $C^\infty(\S_2,\L_2)$.
It follows from Lemma \ref{L:propp} that if $F_i \in C^\infty(\S,\L)$ agrees with $G_i \in C^\infty(\S^2,\L^2)$  on $C^\infty(\S_2,\L_2)$ for $i =1,2$, then  $\{F_1,F_2\}_\S$ agrees with $\{G_1,G_2\}_{\S^2}$ on $C^\infty(\S_2,\L_2)$ as well.
We will call this property of comultiplication {\it Property P} and use it in the definition of a formal symplectic groupoid.
The mappings $s,t: \S \to M$ induce the algebra morphisms
\[
S,T: C^\infty(M) \to C^\infty(\S,\L). 
\]
The source mapping $S$ is a Poisson morphism and the target mapping $T$ is an anti-Poisson morphism. For any $f,g \in C^\infty(M)$ the elements $Sf, Tg \in C^\infty(\S,\L)$ Poisson commute. The unit mapping $\epsilon: M \to \S$ identifies $M$ with $\L$ and thus induces the algebra morphism $E: C^\infty(\S,\L) \to C^\infty(M)$.
Axioms (A3) and (A4) imply that 
\begin{equation}\label{E:ax34}
ES = \id_{C^\infty(M)} \mbox{ and } ET =  \id_{C^\infty(M)}.
\end{equation}
The inverse mapping $i: \S \to \S$ leaves fixed the elements of $\L$ and therefore induces the antisymplectic involutive algebra morphism $I : C^\infty(\S,\L) \to C^\infty(\S,\L)$. It follows from Axiom (A7) that 
\begin{equation}\label{E:ax7}
IS = T.
\end{equation}
To find the formal analogue of multiplication in a symplectic groupoid, we
need a different description of the algebra $C^\infty(\S_n,\L_n)$. For $f \in C^\infty(M)$ introduce functions $S_n^k f, T_n^k f \in C^\infty(\S^n,\L^n)$ by the following formulas:
\begin{equation}\label{E:snktnk}
S_n^k f =  \underset{n}{\underbrace{1 \otimes \ldots \otimes \overset{k-\mbox{th}}{\overbrace{Sf}} \otimes \ldots \otimes1}}, \
T_n^k f =  \underset{n}{\underbrace{1 \otimes \ldots \otimes \overset{k-\mbox{th}}{\overbrace{Tf}} \otimes \ldots \otimes1}}.
\end{equation}
Denote by $\I_n$ the ideal in $C^\infty(\S^n,\L^n)$ generated by the functions $S_n^{k+1}f - T_n^k f, \ f \in C^\infty(M), \ 1 \leq k \leq n-1$.
Taking into account that $\S_n \cap \L^n = \L_n$, we see that the inclusion of $\S_n$ into $\S^n$ induces the following exact sequence of algebras:
\[
    0 \to \I_n \to C^\infty(\S^n,\L^n) \to C^\infty(\S_n,\L_n) \to 0,
\] 
whence $C^\infty(\S_n,\L_n)$ is canonically isomorphic to the quotient algebra $C^\infty(\S^n,\L^n)/\I_n$. Denote
\begin{equation}\label{E:en}
           \E_n := C^\infty(\S^n,\L^n)/\I_n.
\end{equation}
Notice that $\E_1 = C^\infty(\S,\L)$.

{\it Remark.} A formal neighborhood $(X,Y)$ is the simplest example of a formal manifold which, in general, should be defined as a ringed space on $Y$.  If a formal symplectic groupoid is defined as a formal neighborhood $(\S,\L)$, it is too restrictive to require the existence of the manifolds $\S_n$ for $n \geq 2$. This is why from now on we will automatically replace the algebra $C^\infty(\S_n,\L_n)$ with $\E_n$ for $n \geq 2$. In particular, we consider the comultiplication mapping $m^*$ as a mapping from $C^\infty(\S,\L)$ to $\E_2$ and the algebra morphism $\iota_n^*$ as the quotient mapping from $C^\infty(\S^n,\L^n)$ to $\E_n$.

Axioms (A1) and (A2) imply the following identities in the algebra $\E_2$: for $f \in C^\infty(M)$
\begin{equation}\label{E:ax1}
    m^*(Sf) = \iota^*(Sf \otimes 1)
\end{equation}
and
\begin{equation}\label{E:ax2}
    m^*(Tf) = \iota^*(1 \otimes Tf),
\end{equation}
respectively.

In order to state the formal analogues of axioms (A5),(A6),(A8), and (A9) we need one more mapping. Denote by $\delta: \S \to \S \times \S$ the diagonal inclusion of $\S$. Since $\delta(\L)$ is the diagonal of $\L\times \L$, the mapping $\delta$ induces the dual morphism
\[
     \delta^*: C^\infty(\S^2,\L^2) \to C^\infty(\S,\L).
\]
In what follows $F \in C^\infty(\S,\L)$ and $G \in C^\infty(\S^2,\L^2)$ agree on $\E_2$, i.e.,
 $m^* F = \iota^* G$ in $\E_2$. 
Axiom (A5) implies that $F(\alpha) = G(\epsilon(s(\alpha)), \alpha)$, whence 
\begin{equation}\label{E:ax5}
    F = (\delta^* \circ (SE \otimes 1))G.
\end{equation}
Similarly, it follows from Axiom (A6) that
\begin{equation}\label{E:ax6}
    F = (\delta^* \circ (1 \otimes TE))G.
\end{equation}
Axioms (A8) and (A9) imply that
\begin{equation}\label{E:ax89}
    (SE)F = (\delta^* \circ (1 \otimes I))G \mbox{ and } (TE)F = (\delta^* \circ (I \otimes 1))G,
\end{equation}
respectively.

Now we need to state the formal analogue of the associativity of the groupoid multiplication. The mapping $m^*\otimes 1$ maps $F \in C^\infty(\S^2,\L^2)$ to a coset in $C^\infty(\S^3,\L^3)$ of the 
the ideal  generated by  the functions $Tf \otimes 1 \otimes 1 - 1 \otimes Sf \otimes 1$. This ideal belongs to the ideal $\I_3$. Therefore the image of $F$ with respect to the mapping $m^*\otimes 1$ is a well defined element of $\E_3$. It can be checked using formula (\ref{E:ax2}) that  the homomorphism  $m^*\otimes 1$ maps the ideal $\I_2$ to $\I_3$.
This implies that the mapping $m^*\otimes 1$ induces a well defined mapping from $\E_2$ to $\E_3$, which we will denote $(m_2^1)^*$. Similarly, we construct the mapping $(m_2^2)^*:\E_2 \to \E_3$ induced by $1 \otimes m^*$. The associativity of the groupoid multiplication implies that
\begin{equation}\label{E:assoc}
      \left(m_2^1 \right)^* \circ m^* =  \left(m_2^2 \right)^* \circ m^*.
\end{equation}

To define a formal symplectic groupoid over a Poisson manifold $M$, we begin with a collection of the following data: a symplectic manifold  $\S$,  a Lagrangian manifold  $\L \subset \S$,  an embedding  $\epsilon: M \to \S$ such that $\epsilon(M) = \L$, its dual $E: C^\infty(\S,\L) \to C^\infty(M)$, a Poisson morphism $S: C^\infty(M) \to C^\infty(\S,\L)$ and  an anti-Poisson morphism $T: C^\infty(M) \to C^\infty(\S,\L)$ such that $Sf$ and $Tg$ Poisson commute for any $f,g \in C^\infty(M)$, and an involutive antisymplectic automorphism $I : C^\infty(\S,\L) \to C^\infty(\S,\L)$. For $f \in C^\infty(M)$ introduce the functions $S_n^k f, T_n^k f \in C^\infty(\S^n,\L^n)$ by formulas (\ref{E:snktnk}). For each $n$ define the ideal $\I_n$ in $C^\infty(\S^n,\L^n)$ generated by the functions $S_n^{k+1}f - T_n^k f$, where $f \in C^\infty(M)$ and $1 \leq k \leq n-1$,
and the quotient algebra $\E_n = C^\infty(\S^n,\L^n)/\I_n$ as above. Denote by $\iota_n^*: C^\infty(\S^n,\L^n) \to \E_n$ the quotient mapping. There should exist a comultiplication mapping $m^* : C^\infty(\S,\L) \to \E_2$ which has Property P and satisfies the formal analogues of axioms (A1) - (A9) given by formulas (\ref{E:ax1}), (\ref{E:ax2}), (\ref{E:ax34}), (\ref{E:ax5}), (\ref{E:ax6}), (\ref{E:ax7}), and (\ref{E:ax89}), respectively. It should generate the mappings 
\[
 \left(m_2^1\right)^*,  \left(m_2^2\right)^*: \E_2 \to \E_3
\]
as above so that the coassociativity condition (\ref{E:assoc}) is satisfied. In what follows we will refer to the formal analogues of axioms (A1) - (A9) as to axioms (FA1) - (FA9).

\section{Formal symplectic realization of a Poisson manifold}\label{S:formsymp}

If $\S$ is a symplectic manifold, $M$ a Poisson manifold, $s: \S \to M$ a surjective submersion which is a Poisson mapping, and $\epsilon: M\to \S$ an embedding such that $s \circ \epsilon = \id_M$ and $\L = \epsilon(M) \subset \S$ is a Lagrangian manifold, then $\S$ is called a strict symplectic realization of the Poisson manifold $M$ (see \cite{CDW}).
It is known that, given a  strict symplectic realization $\S$ of the Poisson manifold $M$,  there exists a canonical local symplectic groupoid over the manifold $M$ defined on a neighborhood of $\L$ in $\S$, such that $s$ is its source mapping (\cite{CDW}, Thm. 1.2 on page 44).

In this section we will prove a formal version of this theorem. Let $\S$ be a symplectic manifold, $M$ a Poisson manifold, and $\epsilon: M\to \S$ an embedding such that $s \circ \epsilon = \id_M$ and $\L = \epsilon(M) \subset \S$ is a Lagrangian manifold. Denote by $E: C^\infty(\S,\L) \to C^\infty(M)$ the dual mapping of $\epsilon$. Then, if there is given a formal Poisson morphism $S: C^\infty(M) \to C^\infty(\S,\L)$ such that $ES = \id_{C^\infty(M)}$, we say that the formal neighborhood $(\S,\L)$  is a formal strict symplectic realization of the Poisson manifold $M$. 
\begin{theorem}\label{T:main}
  Given a formal strict symplectic realization of a Poisson manifold $M$ on the formal neighborhood $(\S,\L)$ of a Lagrangian submanifold $\L$ of a symplectic manifold $\S$ via a formal Poisson morphism $S: C^\infty(M) \to C^\infty(\S,\L)$, there exists a unique formal symplectic groupoid on $(\S,\L)$ over the manifold $M$ such that $S$ is its source mapping.
\end{theorem}

Assume there is a formal strict symplectic realization $(\S,\L)$ of the Poisson manifold $M$. Given an element $F \in C^\infty(\S,\L)$, denote by $H_F = \{F,\cdot\}_\S$ the formal Hamiltonian vector field corresponding to the formal Hamiltonian $F$. Denote by $\lambda$ the representation of the Lie algebra  $\g := (C^\infty(M), \{\cdot,\cdot\}_M)$ on the space $C^\infty(\S,\L)$ given by the formula
\[
     \lambda(f) = H_{Sf}, 
\]
where $f \in \g$. Extend the representation $\lambda$ to the universal enveloping algebra $\U(\g)$ of $\g$ and define a mapping
\[
    \langle F \rangle : \U(\g) \to C^\infty(M)
\]
by the formula 
\[
    \langle F\rangle (u) = E(\lambda(u)F),
\]
where $u \in \U(\g)$. Denote the multiplication in the algebra $\U(\g)$ by $\bullet$, so that $f \bullet g - g \bullet f = \{f,g\}_M$ for $f,g \in \g$.
We will often work in the following local framework on $\S$. Let $U$ be a Darboux chart on $\S$ with the local coordinates $\{x^k, \xi_l\}$ such that $\L \cap U$ is given by the equations $\xi = 0$ and for $F,G \in C^\infty(U)$
\begin{equation}\label{E:poiss}
     \{F,G\}_\S = \p^k F \,\p_k G - \p^k G\, \p_k F,
\end{equation}
where $\p_k = \p/\p x^k$ and $\p^l = \p/\p \xi_l$. We will say that these Darboux coordinates are {\it standard}. In this framework a  formal function $F \in C^\infty(\S,\L)$ will be represented on the formal neighborhood $(U, \L \cap U)$ as an element of $C^\infty(\L \cap U)[[\xi]]$ and the coordinates $\xi_l$ will be treated as formal variables. We identify $M$ and $\L$ via the mapping $\epsilon$, so that $\{x^k\}$ are used also as coordinates on $\epsilon^{-1}(\L \cap U)$. In particular, for $F=F(x,\xi)$ we have $E(F)(x) = F(x,0)$. A local expression for the Poisson bracket on $M$ is 
\begin{equation}\label{E:poissm}
   \{f,g\}_M = \eta^{ij}\p_i f \, \p_j g,
\end{equation}
where $f,g \in C^\infty(M)$.
\begin{lemma}\label{L:sf}
    Given a function $f \in C^\infty(M)$, the element $Sf \in C^\infty(\S,\L)$ can be written in standard local coordinates $(x,\xi)$ on a Darboux chart $U \subset \S$ as
\begin{equation}\label{E:sfx}
   Sf(x,\xi) = f(x) + \alpha^{ij}(x)\, \p_i f\,  \xi_j \pmod{\xi^2}
\end{equation}
for some function $\alpha^{ij}(x)$ such that $\alpha^{ij} - \alpha^{ji} = \eta^{ij}$.
\end{lemma}
\begin{proof}
Denote $s^i = Sx^i$. Since $E(s^i) = x^i$, expanding $s^i(x,\xi)$ with respect to the formal variables $\xi$ we get that $s^i = x^i + \alpha^{ij}(x) \xi_j \pmod{\xi^2}$ for some function $\alpha^{ij}(x)$. It follows from the fact that $S$ is an algebra morphism, that 
\[
Sf(x,\xi) = f(s(x,\xi)) = f(x) + \alpha^{ij}(x)\, \p_i f\, \xi_j \pmod{\xi^2}. 
\]
Notice that the `substitution'  $f(s(x,\xi))$ is understood as a composition of formal series. Since $S$ is a Poisson morphism, we have that $\{Sf,Sg\}_\S = S(\{f,g\}_M)$ for any $f,g \in C^\infty(M)$. On the one hand, 
according to formulas (\ref{E:poiss}) and (\ref{E:sfx}),
\[
\p^k (Sf) \, \p_k (Sg) - \p^k (Sf) \, \p_k (Sg) = \alpha^{ik} \, \p_i f \, \p_k g - \alpha^{ik} \,  \p_i g\,  \p_k f \pmod{\xi}. 
\]
On the other hand, $S(\{f,g\}_M) = \eta^{ij}\, \p_i f \,  \p_j g \pmod{\xi}$, which concludes the proof.
\end{proof}
\begin{lemma}\label{L:null}
   For any $F \in  C^\infty(\S,\L)$ and $u \in \U(\g)$ the mapping $C^\infty(M) \ni f \mapsto \langle F \rangle (f \bullet u)$ is a derivation on $C^\infty(M)$.
\end{lemma}
\begin{proof}
Let us show that the mapping $C^\infty(M) \ni f \mapsto E(H_{Sf}F)$ is a derivation. Using Lemma \ref{L:sf} and formula (\ref{E:sfx}), we obtain that in local Darboux coordinates
\begin{align}\label{E:loc}
E(H_{Sf}F) = E(\{Sf, F\}_\S) = 
E(\p^k(Sf) \p_k F - \p^k F \p_k (Sf)) =\nonumber \\
 \alpha^{ik}\p_i f \, E(\p_k F) - \p_k f  \, E(\p^k F). 
\end{align}
To prove the statement of the Lemma, the element $F$ should be replaced with $\lambda(u)F$.
\end{proof}

Denote by $\c$ the space of linear mappings $C: \U(\g) \to C^\infty(M)$ such that for any $u \in \U(\g)$ the mapping  $C^\infty(M) \ni f \mapsto C(f \bullet u)$ is a derivation on $C^\infty(M)$. Lemma \ref{L:null} implies that the mapping 
\[
 \chi: F \mapsto \langle F \rangle
\]
maps $C^\infty(\S,\L)$ to $\c$. We will prove that the mapping $\chi : C^\infty(\S,\L) \to \c$ is actually a bijection. Each element $C\in \c$ is completely determined by the family of polydifferential operators $\{C_n\}, \ n \geq 0,$ on $M$, where $C_n$ is the $n$-differential operator such that
\begin{equation}\label{E:corresp}
    C_n(f_1, \ldots, f_n) = C(f_1 \bullet \ldots \bullet f_n).
\end{equation}
The operators $\{C_n\}$ enjoy the following two properties.

\noindent{\it Property A}. Each operator $C_n, n \geq 0,$ is a derivation in the first argument.

\noindent{\it Property B}.  For any $k,n$ such that $1 \leq k \leq n-1$
\begin{align}
    C_n(f_1, \ldots, f_k,f_{k+1}, \ldots f_n) - C_n(f_1, \ldots, f_{k + 1},f_k, \ldots f_n) =\nonumber\\
 C_{n-1}(f_1, \ldots, \{f_k,f_{k+1}\}_M, \ldots f_n).
\end{align}

We will call a family $\{C_n(f_1, \ldots, f_n)\}, n \geq 0$, of polydifferential operators on $M$ {\it coherent} if it has Properties A and B. The correspondence $C \mapsto \{C_n\}$ given by formula (\ref{E:corresp}) is a bijection between the space $\c$ and the set of all coherent families.

It is easy to show that each operator $C_n$ from a coherent family
annihilates constants (i.e., $C_n(f_1,\ldots, f_n) =0$ if $f_k = 1$ for at least one index $k$) and is of order not greater than $k$ in the $k$th argument for $1 \leq k \leq n$.

It is important to notice that if $\{C_n\}, n \geq 0$, is a coherent family on $M$ and $\phi \in C^\infty(M)$, then the operators $\{\phi \cdot C_n\}, n \geq 0$, also form a coherent family. This observation means that one can apply partition of unity arguments to the coherent families. 

The standard increasing filtration on the universal enveloping algebra
$\U(\g)$ induces the dual decreasing filtration $\{\c^{(n)}\}$ on $\c$, i.e., $\c^{(n)}$ consists of all operators $C$ such that the corresponding coherent family $\{C_k\}$ satisfies the condition $C_k = 0$ for $0 \leq k \leq n-1$.  The following lemma is an immediate consequence of Properties A and B of the coherent families.
\begin{lemma}\label{L:calcn}
If $C\in \c^{(n)}$, then $C_n(f_1, \ldots, f_n)$ is a symmetric multiderivation on $M$ (i.e., of order one in each argument and null on constants).
\end{lemma}

We will also consider finite coherent families $\{C_k\},  0 \leq k \leq n$.
It turns out that any $n$-element coherent family can be extended to an $(n+1)$-element coherent family. 
\begin{theorem}\label{T:ext}
Any $n$-element coherent family  $\{C_k\},  0 \leq k \leq n -1$, can be extended to an $(n+1)$-element coherent family $\{C_k\},  0 \leq k \leq n$. The operator $C_n$ is unique up to an arbitrary symmetric multiderivation.
\end{theorem}
We will prove Theorem \ref{T:ext} in the Appendix.

Given an element $F \in C^\infty(\S,\L)$, set $C = \langle F \rangle = \chi(F)$. We will denote the corresponding operator $C_n$ by $\langle F \rangle_n$, so that 
\[
     \langle F \rangle_n (f_1, \ldots, f_n) = E(H_{Sf_1}\ldots H_{Sf_n}F),
\]
where $f_i \in C^\infty(M)$. Denote by $\J = I_\L/I^\infty_\L$ the kernel of the mapping $E: C^\infty(\S,\L) \to C^\infty(M)$, i.e., the ideal of formal functions on $(\S,\L)$ vanishing on $\L$. The powers of this ideal, $\{\J^n\}$, form a decreasing filtration on the algebra $C^\infty(\S,\L)$. Consider an element $F\in \J^n$. The operator $\langle F \rangle_n$ vanishes if $k < n$, therefore $\langle F \rangle \in \c^{(n)}$. Thus the mapping $\chi: C^\infty(\S,\L) \to \c$ is a morphism of  filtered spaces. Notice that the filtrations on  $C^\infty(\S,\L)$ and $\c$ are complete and separated.
\begin{lemma}\label{L:filtr}
Let $F$ be an arbitrary element in $\J^n$. Then
\[
 \langle F \rangle_n (f_1, \ldots, f_n) = E(H_{Sf_1}\ldots H_{Sf_n}F),
\]
$f_i \in C^\infty(M)$, is a symmetric multiderivation which does not depend on the choice of the source mapping $S$. The mapping $\chi: F \mapsto \langle F \rangle_n$ induces an isomorphism of $\J^n/\J^{n+1}$ onto the space of symmetric $n$-derivations on $M$.
\end{lemma}
\begin{proof}
  In standard local Darboux coordinates $(x^k,\xi_l)$ on $\S$ a function $F\in \J^n$ can be written as $F(x,\xi) = F^{i_1\ldots i_n} (x) \xi_{i_1}\ldots \xi_{i_n} \pmod{\xi^{n+1}}$, where $F^{i_1\ldots i_n} (x)$ is symmetric in $i_1, \ldots, i_n$. Taking into account formula (\ref{E:poiss}) and that $Sf = f \pmod{\xi}$, we get that 
\[
 \langle F \rangle_n (f_1, \ldots, f_n) = E(H_{Sf_1}\ldots H_{Sf_n}F) = (-1)^n n! \,F^{i_1\ldots i_n} \p_{i_1}f_1 \ldots \p_{i_n}f_n.
\]
This calculation shows that $F^{i_1\ldots i_n}(x)$ is a symmetric tensor which does not depend on the choice of the source mapping $S$  and that the mapping $\chi: F \mapsto \langle F \rangle_n$ induces an isomorphism of $\J^n/\J^{n+1}$ onto the space of symmetric $n$-derivations on $M$.
\end{proof}

{\it Remark.} One can describe the tensor $F^{i_1\ldots i_n}(x)$ (and the corresponding multiderivation) independently, regardless the existence of the source mapping $S$. The description is based upon the identification of the conormal bundle of $\L \subset \S$ with its tangent bundle $T\L$.

\begin{proposition}\label{P:iso}
The mapping $\chi:C^\infty(\S,\L) \to \c$ is a bijection.
\end{proposition}
\begin{proof}
The mapping $\chi$ is a morphism of complete Hausdorff filtered spaces. According to Theorem \ref{T:ext}, the quotient space $\c^{(n)}/\c^{(n+1)}$ is isomorphic to the space of symmetric $n$-derivations on $M$, the isomorphism being induced by the mapping $\c^{(n)} \ni C \mapsto C_n$. Lemma \ref{L:filtr} thus shows that the mapping $\chi$ induces an isomorphism of $\J^n/\J^{n+1}$ onto $\c^{(n)}/\c^{(n+1)}$, whence the Proposition follows.
\end{proof}

Using Proposition \ref{P:iso} we will transfer the structure of Poisson algebra from $C^\infty(\S,\L)$ to $\c$ via the mapping  $\chi$. It turns out that the resulting Poisson algebra structure on $\c$ does not depend on the mapping $S$ and can be described canonically and intrinsically in terms of the Poisson structure on $M$.

Denote by $\delta_\U : \U(\g) \to \U(\g) \otimes \U(\g)$ the standard cocommutative coproduct on $\U(\g)$, so that $\delta_\U(f) = f \otimes 1 +1 \otimes f$ for $f \in \g$. 

For $F,G \in C^\infty(\S,\L),\ u \in \U(\g),$ we have
\begin{align}\label{E:conv}
   \langle FG \rangle(u) = E(\lambda(u)(FG)) = 
\sum_i E((\lambda(u'_i)F)(\lambda(u''_i)G)) = \nonumber \\
\sum_i E(\lambda(u'_i)F)E(\lambda(u''_i)G) = \sum_i  \langle F \rangle(u'_i) \langle G \rangle (u''_i). 
\end{align}
Here, as well as in the rest of the paper, we use the notation 
\[
\delta_\U (u) = \sum_i  u'_i \otimes u''_i.
\]
For $A,B \in \c$ denote by $AB$ their convolution product, so that
\begin{equation}\label{E:convprod}
    (AB)(u) = \sum_i A(u'_i)B(u''_i).
\end{equation}
This product is commutative since $\delta_\U$ is cocommutative.
Formula (\ref{E:conv}) shows that the mapping $\chi$ is an algebra isomorphism from $C^\infty(\S,\L)$ to $\c$ endowed with the convolution product. 

For $F,G \in C^\infty(\S,\L)$ we obtain by setting $f = E(G)$ in formula (\ref{E:loc}) that
\begin{equation}\label{E:hatf}
   E\left(H_{(SE)(G)}F\right) = \alpha^{ik} E(\p_i G)E(\p_k F)  - E(\p^k F)E(\p_k G).
\end{equation}
Swapping $F$ and $G$ in (\ref{E:hatf}) and subtracting the resulting equation from (\ref{E:hatf}) we get, taking into account formulas (\ref{E:poiss}), (\ref{E:poissm}), and Lemma \ref{L:sf}, that
\begin{align}\label{E:poisse}
   E(\{F,G\}_\S) =  E\left(H_{(SE)(F)}G\right) - E\left(H_{(SE)(G)}F\right) - \nonumber\\\{E(F),E(G)\}_M.
\end{align}
For $u,v \in \U(\g)$ 
\begin{align}\label{E:ehse}
   E\left(H_{(SE)(\lambda(u)F)}\lambda(v)G\right) =  
   E\left(H_{S(\langle F \rangle (u))}\lambda(v)G\right) =
   E\left(\lambda(\langle F \rangle (u))\lambda(v)G\right) =  \nonumber \\
   E\left(\lambda(\langle F \rangle (u) \bullet v)G\right) =\langle G \rangle ( \langle F \rangle (u) \bullet v).
\end{align}
Using the Jacobi identity, formulas (\ref{E:poisse}) and (\ref{E:ehse}),  we obtain that
\begin{align}\label{E:poisshat}
   \langle\{F,G\}_\S\rangle(u) = E\left(\lambda(u) \{F,G\}_\S\right) = \nonumber\\
\sum_i 
E\left(\{\lambda(u'_i)F,\lambda(u''_i)G\}_\S\right) =  \sum_i \big( \langle G \rangle ( \langle F \rangle (u'_i) \bullet u''_i) - \\
 \langle F \rangle ( \langle G\rangle (u''_i) \bullet u'_i) - \{\langle F \rangle (u'_i),\langle G \rangle (u''_i)\}_M \big). \nonumber
\end{align}
Notice that in formula (\ref{E:poisshat}) the functions  $\langle F \rangle (u'_i),\langle G \rangle (u''_i) \in C^\infty(M)$ are used as elements of the Lie algebra $\g$.
Formula (\ref{E:poisshat}) shows that the mapping $\chi$ transfers the Poisson bracket from $C^\infty(\S,\L)$ to the following Poisson bracket on $\c$:
\begin{align}\label{E:poissc}
   \{A,B\}_\c(u) =  \sum_i \Big( B( A(u'_i) \bullet u''_i) - A( B(u''_i) \bullet u'_i) - \nonumber \\
\{A(u'_i), B(u''_i)\}_M \Big),
\end{align}
where $A,B \in \c$. We see that the bracket (\ref{E:poissc}) is defined intrinsically in terms of the Poisson structure on $M$. One can prove that the bracket (\ref{E:poissc}) defines a Poisson algebra structure on the algebra $\c$ regardless the existence of the mapping $S$.

Now we can construct an anti-Poisson morphism $T: C^\infty(M) \to C^\infty(\S,\L)$ such that $ET = \id_{C^\infty(M)}$ and the formal functions $Sf$ and $Tg$ Poisson commute for any $f,g \in C^\infty(M)$. 

Denote by $\epsilon_\U: \U(\g) \to \C$ the counit mapping of the algebra $\U(\g)$, so that $\epsilon_\U({\bf 1}) =1$ and $\epsilon_\U(f) = 0$ for $f \in \g$. Here ${\bf 1}$ is the unity in the algebra $\U(\g)$.
Let ${\bf k}$ denote the trivial representation of the algebra $\U(\g)$ on $C^\infty(M)$, i.e., such that
\[
                      {\bf k}(u)f = \epsilon_\U(u) \cdot f
\]
for $u \in  \U(\g), f \in C^\infty(M)$. For $f \in C^\infty(M)$ consider a mapping $X_f \in \c$ such that
\[
      X_f(u) = {\bf k}(u) f,
\]
where $u \in \U(\g)$. For $f,g \in C^\infty(M)$ and $u \in \U(\g)$ we get from formula (\ref{E:poissc}):
\begin{align*}
  \{X_f,X_g\}_\c (u) = - \left(\sum_i \epsilon_\U(u'_i)\cdot \epsilon_\U(u''_i)\right)\{f,g\}_M = \\
- \epsilon_\U(u) \{f,g\}_M = - X_{\{f,g\}_M}(u).
\end{align*}
Thus the mapping $f \mapsto X_f$ is an anti-Poisson morphism from
$C^\infty(M)$ to $\c$.

Let ${\bf h}$ denote the representation of the Lie algebra $\g$ on $C^\infty(M)$ by the Hamiltonian vector fields, ${\bf h}(f) = \{f,\cdot\}_M, \ f \in \g$. Extend it to $\U(\g)$. It follows from the fact that $ES = \id_{C^\infty(M)}$, that
\begin{equation}\label{E:sfangle}
    \langle Sf \rangle (u) = {\bf h}(u) f.
\end{equation}
For $f,g \in C^\infty(M)$ we get from formula (\ref{E:poissc}) that 
$\langle Sf \rangle$ Poisson commutes with $X_g$:
\begin{align*}
   \{\langle Sf \rangle,X_g\}_\c (u) = \sum_i \left( -\epsilon_\U(u''_i) {\bf h}(g \bullet u'_i)f - \epsilon_\U(u''_i)\{{\bf h}(u'_i)f,g \}_M \right) =\\
 \sum_i \left( -\epsilon_\U(u''_i) {\bf h}(g){\bf h}(u'_i)f + \epsilon_\U(u''_i)\{g, {\bf h}(u'_i)f \}_M \right) = 0.
\end{align*}
Taking into account that the mapping $\chi$ is a Poisson algebra isomorphism of  $C^\infty(\S,\L)$ onto $\c$, we define the mapping $T: C^\infty(M) \to C^\infty(\S,\L)$ as follows. For $f \in C^\infty(M)$ $Tf$  is chosen to be the unique element of $C^\infty(\S,\L)$ such that
\begin{equation}\label{E:tfangle}
   \langle Tf \rangle = X_f.
\end{equation}
We see that the mapping $T: C^\infty(M) \to C^\infty(\S,\L)$ is an anti-Poisson morphism and for any $f,g \in C^\infty(M)$ the formal functions $Sf$ and $Tg$ Poisson commute. Thus the mapping $T$ enjoys the properties of the target mapping. On the other hand, if $T$ is the target mapping of a formal symplectic groupoid on $(\S,\L)$ whose source mapping is $S$, it is straightforward that
\[
      \langle Tf \rangle (u) = {\bf k}(u)f.
\]
It means that the target mapping $T$ is uniquely determined by the source mapping $S$.

In order to construct the inverse mapping $I$ and the comultiplication $m^*$ from the mappings $S$ and $T$, we will consider mappings from tensor powers of $\U(\g)$ to $C^\infty(M)$ which generalize the mappings from $\c$. The space  $\Hom(\U(\g)^{\otimes n},C^\infty(M))$ is endowed with the convolution product defined on its elements $A,B$ as follows:
\begin{align}\label{E:convprodn}
    (AB)(u_1 \otimes \ldots \otimes u_n)  = \\
\sum_{i_1, \ldots, i_n} A(u'_{1i_1} \otimes \ldots \otimes u'_{ni_n})B(u''_{1i_1} \otimes \ldots \otimes u''_{ni_n}), \nonumber
\end{align}
where
\[
           \delta_\U(u_k) = \sum_i u'_{ki} \otimes u''_{ki}.
\]
Denote by $\{\cdot,\cdot\}_{\S^n}$ the Poisson bracket on $C^\infty(\S^n)$ (and on $C^\infty(\S^n,\L^n)$) corresponding to the product Poisson structure. For $F\in C^\infty(\S^n,\L^n)$ let $H_F = \{F, \cdot\}_{\S^n}$ denote the corresponding formal Hamiltonian vector field on $(\S^n, \L^n)$. Introduce representations $\lambda_n^k, \ 0 \leq k \leq n$, of the Lie algebra $\g$ on $C^\infty(\S^n,\L^n)$ by the following formulas:
\[
   \lambda_n^0 (f) = H_{S_n^1 f}, \  \lambda_n^n(f) = - H_{T_n^n f}, \mbox{ and } \lambda_n^k (f) = H_{(S_n^{k+1} f - T_n^k f)}
\]
for $1 \leq k \leq n-1$, where the functions $S_n^k,T_n^k \in C^\infty(\S^n,\L^n)$ are given by formulas (\ref{E:snktnk}). These representations pairwise commute. Notice that in these notations the representation $\lambda$ is denoted $\lambda_1^0$. Denote the representation $\lambda_1^1$ by $\rho$ so that $\rho(f) = - H_{Tf}$ for $f \in \g$. Extend the representations $\lambda_n^k$ to the algebra $\U(\g)$. For $u \in \U(\g)$
\begin{align}\label{E:lambdank}
    \lambda_n^0 (u) = \underset{n}{\underbrace{\lambda(u) \otimes 1 \otimes \ldots \otimes 1}},\  \lambda_n^n (u) = \underset{n}{\underbrace{1 \otimes \ldots \otimes 1 \otimes \rho(u)}},\nonumber\\
\mbox{ and } \quad  \lambda_n^k(u) = \sum_i \underset{n}{\underbrace{1 \otimes \ldots \otimes \overset{k-\mbox{th}}{\overbrace{\rho(u'_i)}} \otimes \lambda(u''_i) \otimes \ldots \otimes 1}},
\end{align}
where $1 \leq k \leq n-1$.
Let $\epsilon_n : M \to \S^n$ denote the composition of the identification mapping from $M$ onto $\L_n$ with the inclusion of $\L_n$ into $\S^n$.
Since $\L_n \subset \L^n$, the mapping $\epsilon_n$ induces the algebra morphism $E_n: C^\infty(\S^n,\L^n)\to C^\infty(M)$. In particular, $\epsilon = \epsilon_1$ and $E = E_1$. After some preparations we will show that the morphism $E_n$ intertwines the representations ${\bf h}$ and $\sum_{k = 0}^n \lambda_n^k$.

We cover the submanifold $\L^n \subset \S^n$ by  Cartesian products $U_1 \times \ldots \times U_n$ of standard Darboux charts $U_i \subset \S$ and use the coordinates $\{x^i_{[k]},\xi_{j[k]}\}$ on the $k$-th factor. In particular, in local coordinates $S_n^k f = (Sf)(x_{[k]}, \xi_{[k]})$ and $T_n^k f = (Tf)(x_{[k]}, \xi_{[k]})$. For a function $F = F(x_{[1]}, \xi_{[1]}, \ldots, x_{[n]},\xi_{[n]})$ on $U_1 \times \ldots \times U_n$ we have $E_n(F) = F(x,0,\ldots, x,0)$. We will use below the following obvious formulas,
\begin{equation}\label{E:easy}
   E_n(f(x_{[k]}) F) = f(x) E_n(F) \mbox{ and } \p_i E_n(F) = \sum_{k = 1}^n E_n(\p_{i[k]}F),
\end{equation}
where $\p_{i[k]} = \p/\p x^i_{[k]}$. It can be proved as in Lemma \ref{L:sf} that in local Darboux coordinates $(x,\xi)$ 
\begin{equation}\label{E:tfx}
   Tf(x,\xi) = f(x) + \alpha^{ji}(x)\, \p_i f\,  \xi_j \pmod{\xi^2},
\end{equation}
where the function $\alpha^{ij}(x)$ is the same as in formula (\ref{E:sfx}).
\begin{lemma}\label{L:leib}
For $f \in C^\infty(M)$ and $F \in C^\infty(\S^n,\L^n)$
\[
   {\bf h}(f)E_n(F) = \sum_{k = 0}^n E_n( \lambda_n^k(f) F).
\]
\end{lemma}
\begin{proof}
Using formulas  (\ref{E:poiss}),(\ref{E:poissm}), (\ref{E:sfx}),(\ref{E:tfx}), and Lemma \ref{L:sf} we get:   
\begin{align*}
 E(H_{(Sf - Tf)}F) = E(\{Sf - Tf, F \}_\S) = \\
E(\p^k (\eta^{ij}\p_i f \xi_j) \p_k F) = {\bf h}(h) E(F).
\end{align*}
Now the Lemma follows from formulas (\ref{E:easy}) and the fact that
\[
  \sum_{k = 0}^n \lambda_n^k(f) =  \sum_{k = 1}^n H_{(S_n^k f - T_n^k f)}.
\]
\end{proof}

Denote by $\c_n$ the subspace of $\Hom(\U(\g)^{\otimes(n+1)},C^\infty(M))$
of the mappings $C$ such that 
\begin{itemize}
\item for any $C \in \c_n,\ u_i \in \U(\g), 0 \leq i \leq n,$ and $k$ satisfying $0 \leq k \leq n$, the mapping
\[
  C^\infty(M) \ni f \mapsto C(u_0 \otimes \ldots \otimes f \bullet u_k \otimes \ldots \otimes u_n)
\]
is a derivation on $C^\infty(M)$; and
\item for any $f \in C^\infty(M)$
\begin{equation}\label{E:propb}
    {\bf h}(f) C(u_0 \otimes  \ldots \otimes  u_n) = \sum_{k = 0}^n C(u_0 \otimes \ldots \otimes  f\bullet u_k \otimes \ldots \otimes  u_n).
\end{equation}
\end{itemize}
The space  $\c_n$ is closed under the convolution product and thus is an algebra. For  an element $F \in C^\infty(\S^n,\L^n)$ define a mapping
\[
 \ll F\rr \in\Hom(\U(\g)^{\otimes(n+1)},C^\infty(M))
\]
such that 
\begin{equation}\label{E:hatfn}
    \ll F \rr (u_0 \otimes  \ldots \otimes  u_n) = E_n(\lambda_n^0(u_0) \ldots \lambda_n^n(u_n)F).   
\end{equation}
A straightforward generalization of the proof of Lemma \ref{L:null} shows that for each $k$ satisfying $0 \leq k \leq n$ the mapping
\[
  C^\infty(M) \ni f \mapsto  \ll F \rr (u_0  \otimes \ldots \otimes  f \bullet u_k \otimes  \ldots \otimes  u_n)
\]
is a derivation on $C^\infty(M)$.  It follows from Lemma \ref{L:leib} that the mapping $C =  \ll F \rr$ satisfies formula (\ref{E:propb}). Thus the mapping
\[
   C^\infty(\S^n,\L^n) \ni F \mapsto \ll F\rr
\]
maps $C^\infty(\S^n,\L^n)$ to $\c_n$. Denote this mapping by $\chi_n$. A simple calculation shows that $\chi_n: C^\infty(\S^n,\L^n) \to \c_n$ is an algebra homomorphism.

Denote by $\tilde \c_n$ the subspace of $\Hom(\U(\g)^{\otimes n},C^\infty(M))$ consisting of the elements
$C \in \tilde \c_n$ such that for any $u_i \in \U(\g), 1 \leq i \leq n,$ and $k$ satisfying $1 \leq k \leq n$, the mapping
\[
  C^\infty(M) \ni f \mapsto C(u_1 \otimes \ldots \otimes  f \bullet u_k \otimes  \ldots  \otimes  u_n)
\]
is a derivation on $C^\infty(M)$. Notice that in these notations $\c = \tilde \c_1$. 
The space $\tilde \c_n$ is also an algebra with respect to the convolution product.

Consider a reduction mapping $C \mapsto \tilde C$ from $\c_n$ to $\tilde \c_n$ defined as follows:
\[
     \tilde C(u_1 \otimes  \ldots \otimes  u_n) = C(u_1 \otimes  \ldots \otimes u_n \otimes {\bf 1}),
\]
where ${\bf 1}$ is the unity in the algebra $\U(\g)$ (which should not be confused with the unit constant $1 \in \g$).
Formula (\ref{E:propb}) implies that the reduction mapping $C \mapsto \tilde C$ is a bijection of $\c_n$ onto $\tilde \c_n$. It is easy to check that the reduction mapping $C \mapsto \tilde C$ is an algebra isomorphism of $\c_n$ onto $\tilde \c_n$.  A straightforward calculation shows that the reduction mapping pulls back the Poisson bracket (\ref{E:poissc}) on $\c = \tilde \c_1$ to the Poisson bracket $\{\cdot,\cdot\}_{\c_1}$ on $\c_1$ defined as follows. For $A,B \in \c_1$ and $u,v \in \U(\g)$
\begin{align}\label{E:poissc1}
\{A,B\}_{\c_1}(u \otimes v) = 
- \sum_{i,j} 
\Big (A\big ((B(u''_i \otimes v''_j)\bullet u'_i) \otimes v'_j\big ) + \nonumber \\
B \big (u''_i \otimes (A(u'_i \otimes v'_j) \bullet v''_j)\big)\Big ),
\end{align}
where
\[
           \delta_\U(u) = \sum_i u'_i \otimes u''_i \mbox{ and }  \delta_\U(v) = \sum_j v'_j \otimes v''_j.
\] 
Recall that in (\ref{E:poissc1}) the functions $A(u'_i \otimes v'_j),B(u''_i \otimes v''_j) \in C^\infty(M)$ are treated as elements of the Lie algebra $\g$.
The right-hand side of formula (\ref{E:poissc1}) is skew-symmetric due to formula (\ref{E:propb}) and cocommutativity of the coproduct $\delta_\U$.

For $F \in C^\infty(\S,\L)$  the mapping $\ll F \rr \in \c_1$ such that
\[
   \ll F\rr  (u \otimes v) = E(\lambda(u)\rho(v)F)
\]
for $u,v \in \U(\g)$ is completely determined by its reduction $\langle F \rangle (u) = E(\lambda(u)F)$. Thus the mapping $\chi_1: C^\infty(\S,\L) \to \c_1$ is a Poisson algebra isomorphism. This isomorphism will be used to introduce the inverse mapping $I$ on $C^\infty(\S,\L)$ in the most transparent way.

A simple calculation shows that for  $f \in C^\infty(M)$ and $u,v \in \U(\g)$
\begin{equation}\label{E:sftf}
    \ll Sf \rr (u \otimes v) = {\bf h}(u){\bf k}(v)f \mbox{ and }  \ll Tf \rr (u \otimes v) = {\bf h}(v){\bf k}(u)f.
\end{equation}
Given a mapping $C : \U(\g) \otimes \U(\g) \to C^\infty(M)$, denote by $C^\dagger$ the mapping from $\U(\g) \otimes \U(\g)$ to $C^\infty(M)$ such that
\[
                C^\dagger(u \otimes v) = C(v \otimes u)
\]
for $u,v \in \U(\g)$. It is easy to check that the mapping $C \mapsto C^\dagger$ leaves invariant the space $\c_1$. Formulas (\ref{E:sftf}) indicate that
\[
          \ll Sf \rr ^\dagger = \ll Tf \rr.
\]
Using formulas (\ref{E:convprodn}) for $n=2$ and (\ref{E:poissc1}) one can readily show that the mapping $C \mapsto C^\dagger$ induces an involutive anti-Poisson automorphism of the Poisson algebra $\c_1$. Define a unique mapping $I$ on $C^\infty(\S,\L)$ such that for $F \in C^\infty(\S,\L)$ and $u,v \in \U(\g)$
\begin{equation}\label{E:inverse}
     \ll I(F) \rr (u \otimes v) = \ll F \rr (v \otimes u).
\end{equation}
It follows that the mapping $I$ is an involutive anti-Poisson automorphism of $C^\infty(\S,\L)$ such that 
\begin{equation}\label{E:istits}
IS = T \mbox{ and } IT = S.
\end{equation}

Now assume that $I$ is the inverse mapping of a formal symplectic groupoid on $(\S,\L)$ over $M$ with the source mapping $S$ (and target mapping $T$). Then for $f \in C^\infty(M)$ and $F \in C^\infty(\S,\L)$
\begin{align*}
  I(\lambda(f)F) = I\left(\{Sf,F\}_\S\right) = - \{ISf, I(F)\}_\S =\\
 -\{Tf,I(F)\}_\S = \rho(f)I(F).
\end{align*}
Therefore $I\circ \lambda(u) = \rho(u) \circ I$ for $u \in \U(\g)$. Since $I$ is involutive, $I \circ \rho(u) = \lambda(u)\circ I$. One can derive from the groupoid axioms that $i \circ \epsilon = \epsilon$. Similarly, for a formal symplectic groupoid, the formula
\[
                             (EI)(F) = E(F),
\]
where $F \in C^\infty(\S,\L)$, holds. Now,
\begin{align*}
  \ll I(F)\rr  (u \otimes v) = E(\lambda(u)\rho(v)I(F)) = E(I(\rho(u)\lambda(v)F)) =\\
E(\rho(u)\lambda(v)F) = E(\lambda(v)\rho(u)F) = \ll F \rr (v \otimes u),
\end{align*}
which means that the inverse mapping $I$ is uniquely determined by the source mapping $S$.

Our next task is to construct the comultiplication of the formal symplectic groupoid from the source and target mappings. Denote by $\I_n$, as above,  the ideal in $C^\infty(\S^n,\L^n)$ generated by the functions 
\begin{equation}\label{E:sminust}
S_n^{k+1}f - T_n^k f, 
\end{equation}
where $f \in C^\infty(M), \ 1 \leq k \leq n-1$, and set $\E_n = 
C^\infty(\S^n,\L^n)/\I_n$ as in formula (\ref{E:en}).
\begin{lemma}\label{L:kernel} 
The representations $\lambda_n^k$ leave invariant the ideal  $\I_n$.  The ideal $\I_n$ is in the kernel of the algebra morphism $\chi_n: C^\infty(\S^n,\L^n)\to \c_n$.
\end{lemma}
\begin{proof}For $f,g \in C^\infty(M)$
\[
\lambda_n^k(f)(S_n^{l+1} - T_n^l) g = \left\{ \begin{array} {cl} (S_n^{k+1} - T_n^k)\{f, g\}_M & \mbox{ if } k = l \\  0 & \mbox{ otherwise, }
\end{array}
\right.
\]
whence we see that the representations $\lambda_n^k$ leave invariant the ideal  $\I_n$. Since $E_n(S_n^k f) =f$ and $E_n(T_n^k f) =f$, we get that $E_n(S_n^{k+1}f - T_n^k f) = 0$. Therefore the ideal $\I_n$ is in the kernel of the algebra morphism $E_n: C^\infty(\S^n,\L^n)\to C^\infty(M)$. 
Now the Lemma follows from formula (\ref{E:hatfn}).
\end{proof}
Lemma \ref{L:kernel} implies that the homomorphism $\chi_n$ factors through $\E_n$. Denote by $\psi_n$ the induced homomorphism from $\E_n$ to $\c_n$. Notice that $\E_1 = C^\infty(\S,\L)$ and $\psi_1 = \chi_1$. It can be obtained by a straightforward generalization of the proof of Proposition \ref{P:iso} that the induced homomorphism $\psi_n: \E_n \to \c_n$ is, actually, an isomorphism. Introduce a mapping 
$\theta : \c_1 \to \c_2$ as follows. For $C \in \c_1$ set
\begin{equation}\label{E:deftheta}
     \theta [C](u \otimes v \otimes w) = {\bf k}(v) C(u \otimes w).
\end{equation}
We define the comultiplication $m^* : \E_1 \to \E_2$ as a pullback of the mapping $\theta$ with respect to the isomorphisms $\psi_1,\psi_2$:
\[
    m^*:= \psi_2^{-1} \circ \theta \circ \psi_1.
\]
Assume that $F \in C^\infty(\S,\L)$ and $G \in C^\infty(\S^2,\L^2)$ agree on $\E_2$, i.e., $m^* F = \iota^* G$ in $\E_2$. This is equivalent to the condition that $\psi_2 (m^* F) = \psi_2(\iota^* G)$ in $\c_2$, where $\iota^*: C^\infty(\S^2,\L^2) \to \E_2$ is the quotient mapping. On the one hand, $\psi_2(\iota^* G) = \chi_2(G)$. On the other hand, $\psi_2 (m^* F) = \theta [\psi_1(F)] = \theta [\chi_1(F)]$. Thus $F$ and $G$ agree on $\E_2$ iff
\begin{equation}\label{E:a5check}
   \ll G \rr (u \otimes v \otimes w) = {\bf k}(v)  \ll F \rr (u \otimes w)
\end{equation}
for any $u,v,w \in \U(\g)$.

Now we will check formula (\ref{E:ax1}), i.e.,  Axiom (FA1). For $f \in C^\infty(M)$ we need to show that $m^*(Sf) = \iota^* (Sf \otimes 1)$ or, equivalently, that
for $u,v,w \in \U(\g)$
\begin{equation}\label{E:a1check}
   \ll Sf \otimes 1\rr (u \otimes v \otimes w) = {\bf k}(v)  \ll Sf \rr (u \otimes w).
\end{equation}
An easy calculation with the use of formulas (\ref{E:lambdank}) and  (\ref{E:sftf}) shows that both sides of (\ref{E:a1check}) equal ${\bf h}(u){\bf k}(v) {\bf k}(w)f$, whence the statement follows.  Formula (\ref{E:ax1}) can be checked similarly.

Axiom (FA3), i.e., the identity $ES = \id_{C^\infty(M)}$, is a part of the definition of a formal strict symplectic realization of the Poisson manifold $M$, and the target mapping $T$ was constructed to satisfy the identity $ET = \id_{C^\infty(M)}$, which is Axiom (FA4).

Our next goal is to check formula (\ref{E:ax5}), i.e., Axiom (FA5). We start with a pair of functions $F\in C^\infty(\S,\L)$ and $G\in C^\infty(\S^2,\L^2)$ which agree on $\E_2$, i.e., satisfy condition (\ref{E:a5check}).

We need to check that formula (\ref{E:ax5}) holds. Applying the isomorphism $\chi$ to the both sides of formula (\ref{E:ax5}), we obtain an equivalent condition:
\begin{equation}\label{E:cond5}
  \ll F \rr (u \otimes w) = E\big(\lambda(u)\rho(v)(\delta^* \circ (SE \otimes 1))G\big).
\end{equation}
It is straightforward that
\[
    (\lambda(u)\rho(v))\circ \delta^* = \sum_{i,j}\delta^*\circ \Big( \big((\lambda(u'_i)\rho(v'_j)\big)\otimes\big((\lambda(u''_i)\rho(v''_j)\big)\Big).
\]
Then, using the fact that
\[
    \lambda(u) \circ S = S \circ {\bf h}(u), \   
    \rho(u) \circ S = S \circ {\bf k}(u),
\]
and Lemma (\ref{L:leib}), we see that 
\[
    (\lambda(u)\rho(v))\circ (SE) = \epsilon_\U(v) \cdot \sum_i (SE) \circ (\lambda(u'_i)\rho(u''_i)).
\]
Finally, taking into account that $E \circ \delta^* = E_2,\ E_2 \circ (SE \otimes 1) = E_2$, and formula (\ref{E:lambdank}), we obtain that 
\begin{align*}
      \ll F \rr (u \otimes v) = \sum_{i,j}{\bf k}(v'_j)E_2\big((\lambda(u'_i)\rho(u''_i))\otimes (\lambda(u'''_i)\rho(v''_j))G\big)= \nonumber \\
\sum_i E_2\big((\lambda(u'_i)\rho(u''_i))\otimes (\lambda(u'''_i)\rho(v))G\big) = \sum_i \ll G \rr (u'_i \otimes u''_i \otimes v),
\end{align*}
where we have used the following notation:
\[
     \big((\delta_\U \otimes 1) \circ \delta_\U\big)(u) = \big((1 \otimes \delta_\U) \circ \delta_\U\big)(u) = \sum_i u'_i \otimes u''_i \otimes u'''_i.
\]
Thus condition (\ref{E:cond5}) is equivalent to the following one:
\begin{equation}\label{E:check5fin}
  \ll F \rr (u \otimes v) = \sum_i \ll  G \rr (u'_i \otimes u''_i \otimes v).
\end{equation}
Formula (\ref{E:check5fin}) is an immediate consequence of (\ref{E:a5check}). 

The remaining axioms of a formal symplectic groupoid can be checked along the same lines.

In order to check Property P of the comultiplication we need the following lemma.
\begin{lemma}\label{L:proppcheck}
   If elements $F \in C^\infty(\S,\L)$ and $G \in C^\infty(\S^2,\L^2)$ agree on $\E_2$, then for any $u,v,w \in \U(\g)$ the elements $\tilde F =\epsilon_\U(v)\cdot (\lambda(u)\rho(v)F)$ and $\tilde G = \lambda_2^0(u)\lambda_2^1(v)\lambda_2^2(w)G$ agree on $\E_2$ as well.
\end{lemma}
\begin{proof}
 We have to show that 
\[
    \ll \tilde G \rr (\tilde u \otimes \tilde v \otimes \tilde w) = {\bf k}(\tilde v)  \ll \tilde F \rr (\tilde u \otimes \tilde w)
\]
for any $\tilde u,\tilde v,\tilde w \in \U(\g)$. It follows immediately from the fact that the representations $\lambda_n^k,\ 0 \leq k \leq n$, pairwise commute.
\end{proof}

Assume that elements $F_i \in C^\infty(\S,\L)$ and $G_i \in C^\infty(\S^2,\L^2)$ agree on $\E_2$ for $i = 1,2$. To check Property P we need to prove that
\[
     \ll \{G_1,G_2\}_{\S^2}\rr (u \otimes v \otimes w) = {\bf k}(v)\ll \{F_1,F_2\}_\S \rr (u \otimes w).
\]
A straightforward calculation with the use of formulas (\ref{E:sfx}), (\ref{E:tfx}), and (\ref{E:easy}) applied to condition (\ref{E:a5check}) with $F = F_i, G = G_i$, where $i = 1,2,$
shows that
\[
    E\left(\{F_1,F_2\}_\S \right) = E_2\left(\{G_1,G_2\}_{\S^2}\right).
\]
 Then it remains to use the Jacobi identity and Lemma \ref{L:proppcheck}.

In order to check the coassociativity of the comultiplication $m^*$ we consider the mappings 
\[
 \left(m_2^1\right)^*,  \left(m_2^2\right)^*: \E_2 \to \E_3
\]
induced by $m^*\otimes 1$ and $1 \otimes m^*$ as in Section \ref{S:def}.
These mappings are well defined due to Axioms (FA1) and (FA2) given by formulas (\ref{E:ax1}) and (\ref{E:ax2}) respectively. Pushing forward the mappings $(m_2^1)^*$ and $(m_2^2)^*$ via the isomorphisms $\psi_2,\psi_3$ we obtain the mappings $\theta_2^1, \theta_2^2 : \c_2 \to \c_3$ such that
\[
     \theta_2^1 = \psi_3 \circ  \left(m_2^1\right)^* \circ \psi_2^{-1},\ \theta_2^2 = \psi_3 \circ  \left(m_2^2\right)^* \circ \psi_2^{-1}.
\]
These mappings act on an element $C \in \c_2$ as follows:
\begin{align*}
     \theta_2^1[C](u\otimes v \otimes w \otimes z) = {\bf k}(v)C(u\otimes w \otimes z), \\
\theta_2^2[C](u\otimes v \otimes w \otimes z) = {\bf k}(w)C(u\otimes v \otimes z).
\end{align*}
Now, both $\theta_2^1\circ \theta$ and $\theta_2^2\circ \theta$ map 
$B\in\c_1$ to an element $D \in \c_3$ such that
\[
     D(u\otimes v\otimes w\otimes z) = {\bf k}(v){\bf k}(w) B(u\otimes z),
\] 
which implies the coassociativity of the coproduct $m^*$.

Assume that there is given a formal symplectic groupoid on $(\S,\L)$ over the Poisson manifold $M$ with the source mapping $S$ and comultiplication $m^*$.
To conclude the proof of Theorem \ref{T:main} we need to prove the following statements.
\begin{lemma}\label{E:efiseg}
If elements $F\in C^\infty(\S,\L)$ and $G \in C^\infty(\S^2,\L^2)$ agree on $\E_2$, then
$E(F) = E_2(G)$.
\end{lemma}
\begin{proof}
Axiom (FA5) given by (\ref{E:ax5}) and formula (\ref{E:sfx})  imply that 
\[
E(F) = E(\delta^* \circ (SE \otimes 1))G) = E_2((SE \otimes 1))G) = E_2(G).
\]
\end{proof}
\begin{proposition}\label{P:uniqueprod}
The mapping $\psi_2 \circ m^* \circ \psi_1^{-1}$ coincides with the mapping $\theta$, given by formula (\ref{E:deftheta})
\end{proposition}
\begin{proof}
Axiom (FA1) of a formal symplectic groupoid given by formula (\ref{E:ax1}) means that the formal functions $Sf$ and $Sf \otimes 1$ agree for all $f \in C^\infty(M)$. Similarly, Axiom (FA2) given by formula (\ref{E:ax2}) means that $Tf$ agrees with $1 \otimes Tf$. Finally, zero constant $0$ agrees with $1 \otimes Sf - Tf \otimes 1$, since the function $1 \otimes Sf - Tf \otimes 1$ is in the ideal $\I_2$ which is the kernel of the mapping $\iota^*$. Property P implies that if $F\in C^\infty(\S,\L)$ agrees with $G \in C^\infty(\S^2,\L^2)$, then
\[
    m^*\big(\lambda(f)F\big) = \iota^*\big(\lambda_2^0(f)G\big), \ m^*\big(\rho(f)F\big) = \iota^* \big(\lambda_2^2(f)G\big),  \iota^* \big(\lambda_2^1(f)G\big) = 0.
\]
Thus for $u,v,w \in \U(\g)$
\begin{equation}\label{E:intert}
    \epsilon_\U(v)m^*\big(\lambda(u)\rho(w)F\big) = \iota^*\big(\lambda_2^0(u)\lambda_2^1(v)\lambda_2^2(w)G\big).
\end{equation}
Taking into account Lemma \ref{E:efiseg}  we obtain from (\ref{E:intert})  that 
\[
     {\bf k}(v)\ll F \rr (u\otimes w) = \ll G \rr (u\otimes v\otimes w),
\]
whence the Proposition follows.
\end{proof}
Proposition \ref{P:uniqueprod} shows that the comultiplication $m^*$ is uniquely defined by the source mapping $S$. This concludes the proof of Theorem \ref{T:main}.

{\it Remark.} Let $M$ be a symplectic manifold. Denote by $\bar M$ a copy of the manifold $M$ endowed with the opposite symplectic structure and by $M_{diag}$ the diagonal of $M \times \bar M$. It follows from the results obtained in \cite{Deq} that, given a formal symplectic groupoid ${\bf G}$ on $(\S,\L)$ over a symplectic manifold $M$ with the source mapping $S$ and target mapping $T$, then the mapping
\[
           S \otimes T: C^\infty(M \times \bar M,M_{diag}) \to C^\infty(\S,\L)
\]
is a formal symplectic isomorphism. It can be easily checked that the mapping $S \otimes T$ induces an isomorphism of the formal pair symplectic groupoid on $(M \times \bar M,M_{diag})$ over $M$ with the groupoid ${\bf G}$.

\section{Isomorphisms of formal symplectic groupoids}\label{S:isom}

Let $\S$ be a symplectic manifold and $\L$ its Lagrangian submanifold which is a copy of a given Poisson manifold $M$. In this section we will consider the formal symplectic groupoids on the formal neighborhood $(\S,\L)$ over $M$. It is known that there exists a local symplectic groupoid over $M$ defined on a symplectic manifold $\S'$. Its unit space $\L'$ is a copy of $M$. One can find a symplectomorphism of a neighborhood $V$ of $\L$ in $\S$ onto a neighborhood $V'$ of $\L'$ in $\S'$ which identifies $\L$ with $\L'$. One can then transfer the local symplectic groupoid on $V'$ to $V$ and induce a formal symplectic groupoid on $(\S,\L)$ over $M$. We are going to describe the space of all formal symplectic groupoids on $(\S,\L)$ over $M$ as a principal homogeneous space of a certain pronilpotent infinite dimensional Lie group.

Let ${\bf G}$ and ${\bf G'}$ be two formal symplectic groupoids on $(\S,\L)$ over $M$ with the source mappings
\[
   S, S' : C^\infty(M) \to C^\infty(\S,\L),
\]
target mappings $T,T'$, and inverse mappings $I,I'$ respectively. Denote by $\chi,\chi': C^\infty(\S,\L) \to \c$ and by $\chi_1,\chi'_1: C^\infty(\S,\L) \to \c_1$ the corresponding Poisson isomorphisms, as introduced in Section \ref{S:formsymp}.  For $F \in C^\infty(\S,\L)$ we use the notations $\langle F \rangle = \chi (F),\ \langle F \rangle' = \chi'(F)$. There exists a unique Poisson automorphism $Q$ of $C^\infty(\S,\L)$ such that
\[
                                  \chi' = \chi \circ Q.
\]
It follows from formulas (\ref{E:sfangle}) and  (\ref{E:tfangle}) that for $f \in C^\infty(M)$
\[
     \langle Sf \rangle =  \langle S'f \rangle' \mbox{ and }    \langle Tf \rangle =  \langle T'f \rangle', 
\]
whence 
\begin{equation}\label{E:sqsprime}
        S = Q \circ S' \mbox{ and }   T = Q \circ T'.
\end{equation}
The isomorphisms $\chi_1,\chi'_1 : C^\infty(\S,\L) \to \c_1$ push forward the corresponding inverse mappings $I$ and $I'$ of the formal symplectic groupoids ${\bf G,G'}$ to the same mapping $C \mapsto C^\dagger$ on $\c_1$. Therefore 
\[
QI' = IQ.
\]
We want to descibe the structure of the automorphism $Q$. The isomorphisms $\chi,\chi'$ respect the filtrations on $C^\infty(\S,\L)$ and $\c$. Therefore, the automorphism $Q$ respects the filtration on $C^\infty(\S,\L)$, i.e., $Q(\J^n) \subset \J^n, n \geq 0,$ where $\J = I_\L/I_\L^\infty$ is the kernel of the unit mapping $E : C^\infty(\S,\L) \to C^\infty(M)$ and $\J^0 := C^\infty(\S,\L)$. One can prove a stronger statement.
\begin{proposition}\label{P:descent}
  The operator $Q - 1: C^\infty(\S,\L) \to C^\infty(\S,\L)$ increases the filtration degree by one, i.e., $(Q-1)\J^n \subset \J^{n+1}, n\geq 0$.
\end{proposition}
\begin{proof}
For an arbitrary element $G \in \J^n$ set $F = Q(G) \in \J^n$. We have that 
  $\langle F \rangle = \langle G \rangle'$ and $\langle F \rangle_k = \langle G \rangle_k = \langle G \rangle'_k = 0$ for all $k < n$. According to Lemma \ref{L:filtr}, $\langle G \rangle_n = \langle G \rangle'_n$, whence $\langle F \rangle_k = \langle G \rangle_k$ for all $k \leq n$. Therefore $(Q-1)G = F - G \in \J^{n+1}$, which concludes the proof.
\end{proof}

For $G \in C^\infty(\S,\L)$ set $F = Q(G)$. Using that $\chi(F) = \chi'(G)$, it is easy to check that in standard local Darboux coordinates $(x,\xi)$ on $\S$
\[
      E(\p^\alpha F) = \Phi^{\alpha \beta}_\gamma(x) \p_\beta E(\p^\gamma G),
\]
where $\alpha, \beta, \gamma$ are multi-indices (recall that $\p_i = \p/\p x^i$ and $\p^j = \p/\p \xi_j$). We see that locally $Q = \Psi^\alpha_\beta(x,\xi) \p_\alpha \p^\beta$, i.e., $Q$ is a formal differential operator on the formal neighborhood $(\S,\L)$. Proposition \ref{P:descent} implies that the operator 
\[
  H := \log Q = \log \big(1 + (Q-1)\big) = \sum_{n = 1}^\infty \frac{(-1)^{n+1}}{n}(Q-1)^n. 
\]
on $C^\infty(\S,\L)$ is correctly defined via a $\J$-adically convergent series and increases the filtration degree by one. Since $Q$ is a Poisson automorphism of $C^\infty(\S,\L)$, the operator $H$ is a derivation of $C^\infty(\S,\L)$ which respects the Poisson bracket. The operator $H$ is a formal vector field on $(\S,\L)$ locally given by the formula
\begin{equation}\label{E:lochamilt}
              H = a^i(x,\xi)\p_i + b_j(x,\xi)\p^j,
\end{equation}
where $a^i = 0 \pmod{\xi}$ and $b_j = 0 \pmod{\xi^2}$, since $H$ increases the filtration degree by one. We want to show that $H$ is a formal Hamiltonian vector field on $(\S,\L)$.
\begin{lemma}\label{L:hamilt}
 A formal vector field $H$ on $(\S,\L)$ respects the Poisson bracket $\{\cdot,\cdot\}_\S$ and increases by one the filtration degree in $C^\infty(\S,\L)$ if and only if there exists a formal Hamiltonian $F \in \J^2$ such that $H =H_F$. If $H = H_F$ for some formal Hamiltonian $F \in \J^2$, then $F$ is defined uniquely. 
\end{lemma}
\begin{proof}
Assume that $H$ is given in standard Darboux coordinates by formula (\ref{E:lochamilt}).
The condition that $H$ respects the Poisson bracket $\{\cdot,\cdot\}_\S$ can be expressed in local coordinates as follows:
\[
     \p^i a^j = \p^j a^i,\ \p_i b_j = \p_j b_i,\ \p_i a^j = -\p^j b_i,
\]
which is equivalent to the fact that the formal one-form $A = a^i d\xi_i - b_j dx^j$ is closed.
Introduce a grading $|\cdot|$ on the differential forms in the variables $x,\xi$ such that $|x|=0,|dx|=0,|\xi|=1,|d\xi|=1$. The differential $d = \p_i dx^i + \p^j d\xi_j$ respects the grading. Denote by $A_q$ the homogeneous component of degree $q$ of the form $A$. Then 
\begin{equation}\label{E:aalpha}
A_q = a^i_{q-1} d\xi_i - b_{jq} dx^j, 
\end{equation}
where $a^i_q$ and $b_{jq}$ denote the homogeneous components of $a^i$ and $b_j$ of degree $q$, respectively. Since $a^i = 0 \pmod{\xi}$ and $b_j = 0 \pmod{\xi^2}$, we see from formula (\ref{E:aalpha}) that the series $A = \sum A_q$ starts with the term $A_2$. The form $A$ is closed iff each homogeneous component $A_q$ is closed. Using the standard homotopy argument involving the Euler operator $\xi_j \p^j$ related to the grading, we get that if $A_q$ is closed, there exists a unique function $F_q(x,\xi)$ homogeneous of degree $q$ in $\xi$ such that $A_q = dF_q$. Now, $F = F_2 + F_3 + \ldots$  is the unique element of $\J^2$ such that $A = dF$, or, equivalently, such that $H = H_F$.
\end{proof}

It follows from Lemma \ref{L:hamilt} that there exists a unique formal function $F \in \J^2$ such that $Q = \exp H_F$. Now assume that ${\bf G}$ is a formal symplectic groupoid on $(\S,\L)$ over $M$ with the source mapping $S$. 
\begin{lemma}\label{L:unique}
  Let $W$ be a Poisson automorphism of $C^\infty(\S,\L)$ such that $E \circ W = E$ and $W \circ S = S$. Then $W$ is the identity automorphism, $W =1$.
\end{lemma}
\begin{proof}
  Since $W$ is a Poisson automorphism of $C^\infty(\S,\L)$ and $W \circ S = S$, we get for $f \in C^\infty(M)$ and $F \in C^\infty(\S,\L)$ that $W(H_{Sf}F) = W(\{Sf,F\}_\S) = \{WSf,WF\}_\S = \{Sf,WF\}_\S = H_{Sf}W(F)$. Therefore $W \circ \lambda(u) = \lambda(u) \circ W$ for any $u \in \U(\g)$. Taking into account that $E \circ W = E$, we obtain that
\[
   \langle F \rangle (u) = E(\lambda(u)F) = E(W\lambda(u)F) = E(\lambda(u)WF) =\langle W(F) \rangle (u).
\]
Proposition \ref{P:iso} implies that $W = 1$, which concludes the proof.
\end{proof}

Take an arbitrary element $F \in \J^2$. The operator $H_F$ on $C^\infty(\S,\L)$ increases the filtration degree by one, therefore there is a Poisson automorphism $Q = \exp H_F$ of $C^\infty(\S,\L)$ such that $E \circ Q =E$. The mapping $S'$ uniquely determined by the equation $S = Q \circ S'$ is a Poisson morphism from $C^\infty(M)$ to $C^\infty(\S,\L)$ with the property that $ES' = \id_{C^\infty(M)}$. Therefore it determines a unique formal symplectic groupoid ${\bf G'}$ on $(\S,\L)$ over $M$ whose source mapping is $S'$. Take $F' \in \J^2$ and set $Q' = \exp H_{F'}$. Lemma \ref{L:unique} implies that if  $S = Q' \circ S'$, then $Q = Q'$ and $F = F'$.

The automorphism $Q$ such that $S = Q \circ S'$ plays the role of the equivalence morphism of the groupoids ${\bf G}$ and ${\bf G'}$. 

Denote by $\g_\S$ the pronilpotent Lie algebra $(\J^2, \{\cdot,\cdot\}_\S)$ and by $G_\S = \exp \g_\S$ the corresponding pronilpotent Lie group. The results of this Section can be combined in the following theorem.
\begin{theorem}\label{T:formiso}
  The space of formal symplectic groupoids over a Poisson manifold $M$ defined on the formal symplectic neighborhood $(\S,\L)$ of a Lagrangian submanifold $\L$ of a symplectic manifold $\S$ is a principal homogeneous space of the group $G_\S$ of formal symplectic automorphisms of $C^\infty(\S,\L)$.
\end{theorem}

Let ${\bf G}$ be a formal  symplectic groupoid over a Poisson manifold $M$ defined on the formal neighborhood $(\T,Z)$ of the zero section $Z$ of the cotangent bundle $\T$. Denote by $\tau$ the antisymplectic involutive automorphism of $\T$ given by the formula $\tau: (x,\xi) \mapsto (x,-\xi)$, where $\{x^i\}$ are local coordinates on $M$ lifted to $\T$ and $\{\xi_j\}$ the dual fibre coordinates on $\T$. It induces the dual antisymplectic involutive morphism $\tau^* : C^\infty(\T,Z) \to C^\infty(\T,Z)$. Let $S,T,I$ be the source, target, and inverse mappings of the groupoid ${\bf G}$, respectively. Since $T: C^\infty(M) \to C^\infty(\T,Z)$ is an anti-Poisson morphism such that $ET = \id_{C^\infty(M)}$, the mapping 
\begin{equation}\label{E:tildestau}
    \tilde S = \tau^* \circ T 
\end{equation}
is a Poisson morphism from $C^\infty(M)$ to $C^\infty(\T,Z)$ such that $E\tilde S = \id_{C^\infty(M)}$. Therefore there exists a unique formal  symplectic groupoid  ${\bf \tilde G}$ on $(\T,Z)$ over $M$ whose source mapping is $\tilde S$. We call ${\bf \tilde G}$ {\it the dual formal  symplectic groupoid of} ${\bf G}$. Theorem \ref{T:formiso} implies that there exists a unique symplectic automorphism $Q \in G_\S$ such that 
\begin{equation}\label{E:sqtildes}
            S = Q \circ \tilde S.
\end{equation}
The automorphism $Q$ is uniquely represented as $Q = \exp H_F$ for some element $F \in \J^2$. Since $T = IS$, we get from formulas (\ref{E:tildestau}) and (\ref{E:sqtildes}) that
\[
               S = Q \circ \tau^* \circ I \circ S.
\]
Set $W:= Q \circ \tau^* \circ I$. One can check that $E\circ Q = E, \  E\circ I = E$, and $E\circ  \tau^* = E$, whence $E \circ W = E$. Since $W$ is a Poisson automorphism of $C^\infty(\T,Z)$, it follows from Lemma \ref{L:unique} that $Q \circ \tau^* \circ I = W =1$. Taking into account that the inverse mapping  $I$ is involutive, we obtain that
\begin{equation}\label{E:isomf}
             I = Q \circ \tau^* = \exp H_F \circ \tau^*.
\end{equation}
The Hamiltonian $F$ is canonically related to the formal groupoid ${\bf G}$. Since $\tau^*$ is involutive, we get that $Q \circ \tau^* = \tau^* \circ Q^{-1}$, whence $H_F \circ \tau^* = - \tau^* \circ H_F$, which means that
\begin{equation}\label{E:even}
\tau^* F = F, 
\end{equation}
i.e., that $F(x,\xi) = F(x,-\xi)$.

\section{Canonical formal symplectic groupoid of a natural deformation quantization}\label{S:cannat}

Let $(M,\{\cdot,\cdot\}_M)$ be a Poisson manifold. Denote
by $C^\infty(M)[[\nu]]$ the space of formal series in $\nu$ with
coefficients from $C^\infty(M)$. As introduced in \cite{BFFLS}, a formal
differentiable deformation quantization on $M$ is an associative algebra
structure on $C^\infty(M)[[\nu]]$ with the $\nu$-linear and $\nu$-adically
continuous product $\ast$ (named star-product) given on $f,g\in
C^\infty(M)$ by the formula
\begin{equation} \label{E:star} 
f \ast g = \sum_{r = 0}^\infty \nu^r C_r(f,g), 
\end{equation} 
where $C_r,\, r\geq 0,$ are bidifferential operators on $M$, $C_0(f,g) =
fg$ and $C_1(f,g) - C_1(g,f) = \{f,g\}$.  We adopt the convention that
the unity of a star-product is the unit constant. Two differentiable
star-products $\ast,\ast'$ on a Poisson manifold $(M,\{\cdot,\cdot\}_M)$ are called equivalent if there exists an isomorphism of algebras $B: (C^\infty(M)[[\nu]],\ast') \to (C^\infty(M)[[\nu]],\ast)$ of the form $B= 1 + \nu B_1 + \nu^2 B_2 + \dots,$ where $B_r, r\geq 1,$ are differential
operators on $M$. The existence and classification problem for deformation quantization was first solved in the non-degenerate (symplectic) case (see \cite{DWL}, \cite{OMY}, \cite{F1} for existence proofs and \cite{F2}, \cite{NT}, \cite{D}, \cite{BCG}, \cite{X} for classification) and then Kontsevich \cite{K} showed that every Poisson manifold admits a deformation quantization and that the equivalence classes of deformation quantizations can be parameterized by the formal deformations of the Poisson structure. 

All the explicit constructions of star-products enjoy the following property: for all $r \geq 0$ the bidifferential operator $C_r$ in (\ref{E:star}) is of order not greater than $r$ in each argument (most important examples are Fedosov's star-products on symplectic manifolds and  Kontsevich's star-product on $\R^n$ endowed with an arbitrary Poisson bracket). The star-products with this property were called natural by Gutt and Rawnsley in \cite{GR}, where general properties of such star-products were studied.

Let $\D = \D(M)$ be the algebra of differential operators with smooth comlex-valued coefficients  and $\D[[\nu]]$ be the algebra of formal differential operators on $M$. The algebra $\D$ has a natural increasing filtration $\{\D_r\}$, where $\D_r$ is the space of differential operators of order not greater than $r$. We call a formal differential operator $A = A_0 + \nu A_1 + \dots \in \D[[\nu]]$  {\it natural} if  $A_r \in \D_r$ for any $r \geq 0$. The natural formal differential operators form an algebra which we denote by $\N$.

Let $\T$ be the cotangent bundle of the manifold $M$ and $Z$ be its zero section. Denote by $\epsilon: M \to \T$ the composition of the identifying mapping from $M$ onto $Z$ with the inclusion mapping of $Z$ into $\T$. It induces the dual mapping $E: C^\infty(\T,Z) \to C^\infty(M)$.

If $\{x^k\}$ are local coordinates on $M$ and $\{\xi_k\}$ are the dual fibre coordinates on $\T$, then the principal symbol of  an operator $A \in \D_r$, whose leading term is $a^{i_1\ldots i_r}(x)\p_{i_1}\ldots \p_{i_r}$, is given by the formula $\s_r(A) = a^{i_1\ldots i_r}(x)\xi_{i_1}\ldots \xi_{i_r}$. It is globally defined on $\T$ and fibrewise is a  homogeneous polynomial of degree $r$. We define a $\sigma$-symbol of a natural formal differential operator $A = A_0 + i\nu A_1 + (i\nu)^2 A_2 +\dots$ as the formal series $\sigma(A) = \s_0(A_0) + \s_1(A_1) + \dots$. Such a formal series can be treated as a formal function from $C^\infty(\T,Z)$.  The mapping $\sigma: A \mapsto \sigma(A)$ is an algebra morphism from $\N$ to $C^\infty(\T,Z)$. Moreover, for $A,B \in \N$ the operator $\frac{1}{\nu}[A,B]$ is also natural and
\begin{equation}\label{E:commut}
  \sigma\left(\frac{1}{\nu}[A,B]\right) = \{\sigma(A),\sigma(B)\}_\T,
\end{equation}
where $\{\cdot,\cdot\}_\T$ denotes the standard Poisson bracket on $\T$ and the induced bracket on $C^\infty(\T,Z)$ given locally by the formula 
\[
      \{\Phi,\Psi\}_\T = \p^i \Phi \, \p_i \Psi - \p^i \Psi\,  \p_i \Phi.
\]
For $f,g \in C^\infty(M)[[\nu]]$ denote by $L_f$ and $R_g$ the operators of $\ast$-multi\-pli\-cation by $f$ from the left and of $\ast$-multi\-pli\-cation by $g$ from the right respectively, so that $L_f g = f \ast g = R_g$. The associativity of $\ast$ is equivalent to the fact that $[L_f, R_g]=0$. A star-product $\ast$ on $M$ is natural iff for any  $f,g \in C^\infty(M)[[\nu]]$ the operators $L_f, R_g$ are natural. It was proved in \cite{Deq} that the mappings
\[
          S,T: C^\infty(M) \to C^\infty(\T,Z)
\]
defined by the formulas
\[
     Sf = \sigma(L_f), \ Tf = \sigma(R_f),
\]
where $f \in C^\infty(M)$, are a Poisson and an anti-Poisson morphisms, respectively, which satisfy the formulas $ES = \id_{C^\infty(M)}$ and $ET = \id_{C^\infty(M)}$. Moreover, for $f,g \in C^\infty(M)$ the formal functions $Sf,Tg$ Poisson commute.  For each natural deformation quantization on $M$ we constructed in \cite{Inv} an involutive antisymplectic automorphism $I$ of the Poisson algebra $C^\infty(\T,Z)$ such that $IS = T$ and $IT = S$. It follows from Theorem \ref{T:main} that there exists a canonical formal symplectic groupoid on $(\T,Z)$ over $M$ with the source mapping $S$, target mapping $T$, and inverse mapping $I$. We call it {\it the formal symplectic groupoid of the natural deformation quantization}. 

If $*$ and $*'$ are two equivalent natural star products on $M$, it was proved in \cite{GR}
that any equivalence operator $B$ of these star products satisfying the identity
\[
           Bf * Bg = B(f *' g)
\]
can be represented as $B = \exp \frac{1}{\nu}X$, where $X$ is a natural operator such that $X = 0 \pmod{\nu^2}$. Let ${\bf G}$ and  ${\bf G'}$ be the formal symplectic groupoids of the star products $*$ and $*'$ with the source mappings $S$ and $S'$, respectively. It is easy to check that if $Q$ is the equivalence morphism of these groupoids such that $S = Q \circ S'$, then
\[
            Q = \exp H_{\sigma(X)}.
\]

\section{Deformation quantizations with separation of variables}\label{S:defquant}

Let $M$ be a complex manifold endowed with a Poisson tensor $\eta$ of type (1,1) with respect to the complex structure. We call such manifolds K\"ahler-Poisson. If $\eta$ is nondegenerate, $M$ is a K\"ahler manifold.

If $U \subset M$ is a coordinate chart with local holomorphic coordinates $\{z^k, \bar z^l \}$, we will write $\eta = g^{\bar l k} \bar\p_l  \wedge \p_k$ on $U$, where $\p_k = \p/\p z^k$ and $\bar\p_l = \p/\p \bar z^l$. The condition that $\eta$ is a Poisson tensor is expressed in terms of $g^{\bar l k}$ as follows:
\begin{equation}\label{E:kp}
  g^{\bar l k} \p_k g^{\bar n m} = g^{\bar n k} \p_k g^{\bar l m} \mbox{ and } g^{\bar l k} \bar\p_l g^{\bar n m} = g^{\bar l m} \bar\p_l g^{\bar n k}.
\end{equation}
The corresponding Poisson bracket  on $M$ is given locally as
\begin{equation}\label{E:poissmcompldef}
 \{\phi,\psi\}_M = g^{\bar l k} (\bar\p_l \phi\, \p_k \psi - \bar\p_l\, \psi \p_k\phi).
\end{equation}
We say that a star-product (\ref{E:star}) on a K\"ahler-Poisson manifold $M$ defines a deformation quantization with separation of variables on  $M$ if the bidifferential operators $C_r$ differentiate their first argument in antiholomorphic directions and its second argument in holomorphic ones. 

With the assumption that the unit constant 1 is the unity of the star-algebra $(C^\infty(M)[[\nu]], *)$, the condition that $*$ is a star-product with separation of variables can be restated as follows. For any local holomorphic function $a$ and antiholomorphic function $b$ the operators $L_a$ and $R_b$ are the operators of point-wise multiplication by the functions $a$ and $b$ respectively, $L_a = a, \ R_b = b$.  In such a case it is easy to check that $C_1(\phi,\psi) = g^{\bar l k} \bar\p_l \phi\, \p_k \psi$, so that
\begin{equation}\label{E:spsv}
     \phi * \psi = \phi\psi + \nu g^{\bar l k} \bar\p_l \phi\, \p_k \psi + \dots
\end{equation}
Deformation quantizations with separation of variables on a K\"ahler manifold $M$ (also known as deformation quantizations of the Wick type, see \cite{BW}) are completely described and parameterized by the formal deformations of the K\"ahler form on $M$ in \cite{CMP1}. 
If $\left(g^{\bar l k}\right)$ is an arbitrary matrix with constant entries, the formula 
\[
     (\phi * \psi)(z,\bar z) = \left(\exp\, \nu g^{\bar l k} \frac{\p}{\p \bar v^l}\frac{\p}{\p v_k}\right) \phi(z,\bar v) \psi(v,\bar z)|_{v=z,\bar v=\bar z}
\]
defines a star-product with separation of variables on the K\"ahler-Poisson manifold $(\C^d, g^{\bar l k} \bar \p_l \wedge \p_k)$.
One can give more elaborate examples of deformation quantizations with separation of variables on K\"ahler-Poisson manifolds. We conjecture that  star-products with separation of variables exist on an arbitrary K\"ahler-Poisson manifold and they can be parameterized by the formal deformations of the K\"ahler-Poisson tensor $\eta$ (not the equivalence classes, but the star-products themselves). The nature of this parameterization must be very different from that of the parameterization by the formal deformations of the K\"ahler form in the K\"ahler case (see also \cite{Third}).

For a given star-product with separation of variables $*$ on $M$ there exists a unique formal differential operator $B$ on $M$ such that
\begin{equation}\label{E:ber}
B(ab) = b*a 
\end{equation}
for any local holomorphic function $a$ and antiholomorphic function $b$. The operator $B$ is called the formal Berezin transform (see \cite{Tr}). 
One can check that the operator $\Delta$ defined locally by the formula $g^{\bar l k} \p_k\bar \p_l$ is coordinate invariant and thus globally defined on $M$ and that 
\begin{equation}\label{E:berexpan}
B = 1 + \nu \Delta + \dots. 
\end{equation}
In particular, $B$ is invertible. Introduce a dual star product $\tilde *$ on $M$ by the formula
\begin{equation}\label{E:dual}
    \phi \tilde * \psi = B^{-1}(B\psi * B\phi).
\end{equation}
We will show that $\tilde *$ is a deformation quantization with separation of variables on the complex manifold $M$ endowed with the opposite Poisson tensor $-\eta$. This statement was proved in the K\"ahler case in \cite{Tr}, but the proof does not work in the K\"ahler-Poisson case.

It follows from (\ref{E:ber}) that 
\begin{equation}\label{E:baabbb}
Ba = a \mbox{ and } Bb = b.
\end{equation}
In particular, $B1 = 1$. 
\begin{lemma}\label{L:bab}
  For any local holomorphic function $a$ and antiholomorphic function $b$
\[
     BaB^{-1} = R_a \mbox{ and } BbB^{-1} = L_b.
\]
\end{lemma}
\begin{proof}
We need to show that $BaB^{-1}f = f *a$ for any formal function $f$. Since $B$ is invertible, the function $f$ can be representad as $f = Bg$ for some formal function $g$. Now we need to check that $B(ag) = Bg * a$ for an arbitrary formal function $g$. It suffices to check it only for $g$ of the form $g = \tilde ab$, where $\tilde a$ is a local holomorphic function and $b$ a local antiholomorphic function. We have
\[
  B(a\tilde a b) = b * (\tilde a a) = b * (\tilde a * a) = (b * \tilde a) * a = B(\tilde a b) * a.
\] 
The formula $BbB^{-1} = L_b$ can be proved similarly.
\end{proof}
Denote by $\tilde L_\phi$ the operator of star-multiplication by a function $\phi$ from the left and by $\tilde R_\psi$ the operator of star-multiplication by a function $\psi$ from the right with respect to the star-product $\tilde *$. It follows form (\ref{E:dual}) that 
\begin{equation}\label{E:tildelr}
  \tilde L_\phi = B^{-1} R_{B\phi} B  \mbox{ and } \tilde R_\psi = B^{-1} L_{B\psi} B.
\end{equation}
\begin{proposition}\label{P:dual}
The dual star-product $\tilde *$ given by formula (\ref{E:dual}) is a deformation quantization with separation of variables on the manifold $M$ endowed with the same complex structure but with the opposite Poisson tensor $-\eta$.
\end{proposition}
\begin{proof}
Lemma \ref{L:bab} and formulas (\ref{E:baabbb}) and (\ref{E:tildelr}) imply that for any local holomorphic function $a$ 
\[
  \tilde L_a = B^{-1} R_{Ba} B = B^{-1} R_a B = B^{-1} (BaB^{-1}) B = a.
\]
Similarly, $\tilde R_b = b$ for any local antiholomorphic function $b$. Thus $\tilde *$ is a star-product with separation of variables. Using formulas (\ref{E:spsv}), (\ref{E:berexpan}), and (\ref{E:dual}) we get that
\[
   \phi \tilde * \psi = \phi\psi - \nu g^{\bar l k}\bar \p_l \phi \p_k \psi + \dots,
\]
which implies that $\tilde *$ is a star-product on the K\"ahler-Poisson manifold $(M, -\eta)$.
\end{proof}
Lemma \ref{L:bab} and formula (\ref{E:baabbb}) imply that for any local holomorphic functions $a,\tilde a$ and antiholomorphic functions $b,\tilde b$
\begin{equation}\label{E:babcomm}
     [BaB^{-1},\tilde a] = [R_a,L_{\tilde a}] = 0  \mbox{ and } [BbB^{-1},\tilde b] = [L_b,R_{\tilde b}] = 0. 
\end{equation}

It follows from formula (\ref{E:berexpan}) that $B = \exp\big(\frac{1}{\nu}X\big)$ for some formal differential operator 
\begin{equation}\label{E:operx}
  X = \nu^2 X_2 + \nu^3 X_3 + \dots, 
\end{equation}
where $X_2 = \Delta$. We want to show that the operator $X$ is natural. To this end we need the following technical lemma. If $U$ is a holomorphic chart  on $M$ with local coordinates $\{z^k,\bar z^l\}$ we denote by $\{\zeta_k, \bar \zeta_l\}$ the dual fibre coordinates  on $T^*U$ and set $\p^k = \p/\p\zeta_k$ and $\bar \p^l = \p/\p\bar\zeta_l$.
\begin{lemma}\label{L:diffp}
Given an integer $n \geq 2$, let $X$ be a nonzero differential operator on a holomorphic chart $U$ with coordinates $\{z^k,\bar z^l\}$,  such that the operators $[[X,z^i],z^k]$ and $[[X, \bar z^j],\bar z^l]$ are of order not greater than $n-2$ for any $i,j,k,l$. Then the operator $X$ is of order not greater than $n$.
\end{lemma}
\begin{proof}
Assume that $X$ is a differential operator of order $N > n$.  Its principal symbol $p(\zeta,\bar \zeta)$ is a nonzero homogeneous polynomial of degree $N$ with respect to the fibre coordinates $\{\zeta_k, \bar \zeta_l\}$. The condition that the operator  $[[X,z^i],z^k]$ is of order not greater than $n-2$ means that the function $\p^i \p^k p$ is a polynomial of order not greater than $n-2$ in the variables $\zeta,\bar \zeta$. On the other hand, $\p^i\p^k p$ is of order $N - 2 > n-2$ which means that $\p^i \p^k p = 0$ for any $i,k$. Similarly, $\bar \p^j\bar \p^l p=0$ for any $j,l$. 
Since $N \geq 3$,  at least one of the partial derivatives $\p^i \p^k p$ or $\bar \p^j\bar \p^l p$ should be nonzero. Thus the assumption that $N >n$ leads to a contradiction.
\end{proof}

Formula (\ref{E:baabbb}) implies that for any $n$ the operator $X_n$ in (\ref{E:operx}) annihilates holomorphic and antiholomorphic functions. In particular, $X_n1 =0$. We get from formula (\ref{E:babcomm}) that
\begin{align}\label{E:adx}
    \left[\exp \left(\frac{1}{\nu}\ad X\right) a,\tilde a\right] =   [BaB^{-1},\tilde a] = 0 \mbox{ and } \nonumber \\
\left[\exp \left(\frac{1}{\nu}\ad X\right) b,\tilde b\right] = [BbB^{-1},\tilde b] = 0. 
\end{align}
Expanding the left-hand sides of formulas (\ref{E:adx}) in the formal series in the parameter $\nu$ and equating the coefficient at $\nu^{n-1}$ to zero, we get
\begin{equation}\label{E:sumad}
      \sum_{k \geq 1} \frac{1}{k!} \sum_{i_1 + \ldots + i_k - k = n-1} \left[\left[X_{i_1},\ldots, \left[X_{i_k}, a\right] \ldots \right], \tilde a\right] =0
\end{equation}
for $n \geq 2$. Since all the indices $i_j$ in (\ref{E:sumad}) satisfy the condition $i_j \geq 2$, we have that $n - 1 = i_1 + \ldots + i_k - k \geq k$. Thus we obtain from (\ref{E:sumad}) that
\begin{equation}\label{E:xaa}
   \left[\left[X_n,a\right],\tilde a\right] = -  \sum_{k = 2}^{n - 1} \frac{1}{k!} \sum_{i_1 + \ldots + i_k - k = n-1} \left[\left[X_{i_1},\ldots, \left[X_{i_k}, a\right] \ldots \right], \tilde a\right].
\end{equation}
Similarly,
\begin{equation}\label{E:xbb}
   \left[\left[X_n,b\right],\tilde b\right] = -  \sum_{k = 2}^{n - 1} \frac{1}{k!} \sum_{i_1 + \ldots + i_k - k = n-1} \left[\left[X_{i_1},\ldots, \left[X_{i_k}, b\right] \ldots \right], \tilde b\right].
\end{equation}
The right-hand sides of equations (\ref{E:xaa})  and (\ref{E:xbb}) depend only on $X_k$ with $k < n$. We know that $X_2 = \Delta$ is of order (not greater than) two. Assume that we have proved that $X_k$ is of order not greater than $k$ for all $k < n$. We see from (\ref{E:xaa})  and (\ref{E:xbb}) that $ [[X_n,a],\tilde a]$ and $[[X_n,b],\tilde b]$ are of order not greater than $n-2$. It follows from Lemma \ref{L:diffp} that $X_n$ is of order not greater than $n$. The induction shows that $X$ is indeed a natural operator. We have proved the following proposition.
\begin{proposition}\label{P:bertrans}
  The formal Berezin transform $B$ of a deformation quantization with separation of variables on a K\"ahler-Poisson manifold is of the form $B = \exp \frac{1}{\nu}X$, where $X$ is a natural differential operator such that $X = 0 \pmod{\nu^2}$.
\end{proposition}

It follows from Proposition \ref{P:bertrans} that the conjugation of the formal differential operators with respect to the formal Berezin transform, $A \mapsto BAB^{-1}$, leaves invariant the algebra $\N$ of natural differential operators. In particular, the operators $R_a = BaB^{-1}$ and $L_b = BbB^{-1}$ are natural. Now, if $f = ab = a * b$ we see that $L_f = L_{a * b} = L_a L_b = aL_b$ and $R_f = R_{a * b} = R_b R_a = bR_a$ are natural differential operators. Using the same arguments as in Proposition 1 of \cite{Inv} we can prove the following theorem.
\begin{theorem}\label{T:natural}
 Any deformation quantization with separation of variables on a K\"ahler-Poisson manifold is natural.
\end{theorem}
Theorem \ref{T:natural} was proved in \cite{BW} and \cite{N} in the K\"ahler case.

It follows from Theorem \ref{T:natural} that to any deformation quantization with separation of variables on a K\"ahler-Poisson manifold $M$ there corresponds a canonical formal symplectic groupoid on $(\T,Z)$ over $M$. Since for any deformation quantization with separation of variables $L_a = a$ and $R_b = b$, we see that $Sa = \sigma(L_a) = \sigma(a) = a$ and, similarly, $Tb =b$ (abusing notations we denote by $a$ and $b$ both local functions on $M$ and their lifts to $\T$ with respect to the standard bundle projection). 

Given a K\"ahler-Poisson manifold $M$, we call a formal symplectic groupoid on $(\T,Z)$ over $M$ such that $Sa = a$ and $Tb =b$ for any local holomorphic function $a$ and antiholomorphic function $b$, a {\it formal symplectic groupoid with separation of variables}.

\section{Formal symplectic groupoid with separation of variables}\label{S:formgroupsv}

In this section we will show that for any K\"ahler-Poisson manifold $M$ there is a unique formal symplectic groupoid with sepration of variables over $M$.
Let $U \subset M$ be an arbitrary coordinate chart  with local holomorphic coordinates $\{z^k, \bar z^l \}$. Introduce differential operators $D^k, \bar D^l$ on $U$ by the formulas
\[
D^k \psi = g^{\bar l k} \bar \p_l \psi = - \{z^k, \psi \}_M \mbox{ and } \bar D^l \psi = g^{\bar l k} \p_k \psi = \{\bar z^l, \psi\}_M,
\]
where the Poisson bracket $\{\cdot,\cdot\}_M$ is given by formula (\ref{E:poissmcompldef}).
Conditions (\ref{E:kp}) are equivalent to the statement that 
\begin{equation}\label{E:dbard}
[D^k, D^m]=0 \mbox{ and }[\bar D^l, \bar D^n] =0
\end{equation}
for any $k,l,m,n$. Using the operators $D^k, \bar D^l$ we can write
\begin{equation}\label{E:poissmcompl}
     \{\phi,\psi\}_M = D^k \phi \, \p_k \psi - D^k \psi \, \p_k \phi = \bar\p_l \phi\,\bar D^l \psi - \bar\p_l \psi\,\bar D^l \phi.
\end{equation}
Denote by $\{\zeta_k, \bar \zeta_l\}$ the fibre coordinates on $T^*U$ dual to $\{z^k, \bar z^l \}$. The standard Poisson bracket on $\T$ can be written locally as
\begin{equation}\label{E:poisstm}
    \{\Phi, \Psi\}_\T = \p^k \Phi \, \p_k \Psi - \p^k \Psi \, \p_k \Phi + \bar \p^l \Phi \, \bar \p_l \Psi - \bar \p^l \Psi \, \bar \p_l \Phi, 
\end{equation}
where $\p^k = \p/\p \zeta_k, \ \bar \p^l = \p/\p \bar \zeta_l$. The Poisson bracket  on $\T$ induces a Poisson bracket on $C^\infty(\T,Z)$ which will be denoted also by $\{\cdot,\cdot\}_\T$. 
 Introduce mappings $S,T: C^\infty(U) \to C^\infty(T^*U,Z\cap U)$ by the formulas
\begin{equation}\label{E:st}
   (S\phi)(z,\bar z, \zeta) = e^{\zeta_k D^k}\phi, \ (T\psi)(z, \bar z, \bar \zeta) = e^{\bar \zeta_l \bar D^l}\psi,
\end{equation}
where $\phi,\psi \in C^\infty(M)$, the variables $\zeta,\bar \zeta$ are used as formal parameters, and the exponentials are defined via formal Taylor series.
\begin{proposition}\label{P:pantip}
   The mappings 
\[
S,T: (C^\infty(U), \{\cdot,\cdot\}_M) \to (C^\infty(T^*U,Z\cap U),\{\cdot,\cdot\}_\T)
\]
 are a Poisson and an anti-Poisson morphisms, respectively. For any $\phi,\psi \in C^\infty(U)$ the elements $S\phi,T\psi \in C^\infty(T^*U,Z\cap U)$ Poisson commute.
\end{proposition}
\begin{proof} Since $D^k,\bar D^l$ are derivations of the algebra $C^\infty(T^*U,Z\cap U)$, the operators $e^{\zeta_k D^k},e^{\bar \zeta_l\bar D^l}$ are automorphisms of this algebra
which implies that $S,T$ are algebra homomorphisms. We see from (\ref{E:dbard}) and (\ref{E:st}) that
\begin{equation}\label{E:pspt}
     \p^k (S\phi) = D^k (S\phi) \mbox{ and } {\bar\p}^l (T\psi) = {\bar D}^l (T\psi).
\end{equation}
Fix arbitrary functions $\phi,\psi \in C^\infty(U)$ and introduce an element $u(\zeta)\in C^\infty(T^*U, Z\cap U)$ by the formula
\[
         u(\zeta) = \{S\phi,S\psi\}_\T.
\]
In order to show that $S$ is a Poisson morphism we need to prove that $u(\zeta) = S\{\phi,\psi\}_M = e^{\zeta_k D^k}\{\phi,\psi\}_M$. This amounts to checking that $u(0) = \{\phi,\psi\}_M$ and that $\p^m u = D^m u$. Using (\ref{E:poisstm}), (\ref{E:st}), and (\ref{E:pspt}) we get
\begin{align}\label{E:uzeta}
  u(\zeta) = \{S\phi,S\psi\}_\T = \p^k(S\phi)\p_k(S\psi) - \p^k(S\psi)\p_k(S\phi) \nonumber \\
= D^k(S\phi)\p_k(S\psi) - D^k(S\psi)\p_k(S\phi).
\end{align}
It follows from (\ref{E:poissmcompl}) and (\ref{E:uzeta}) that 
\begin{equation}\label{E:uinit}
u(0) = D^k\phi\, \p_k\psi - D^k\psi\,\p_k\phi = \{\phi,\psi\}_M.
\end{equation}
Now, taking into account (\ref{E:kp}) and (\ref{E:dbard}), we obtain from (\ref{E:uzeta}) that
\begin{align*}
   \p^m u - D^m u = \big(D^mD^k(S\phi)\p_k(S\psi) - D^mD^k(S\psi)\p_k(S\phi) + \\
 D^k(S\phi)\p_k(D^mS\psi) - D^k(S\psi)\p_k(D^mS\phi)\big) - \\
\big(D^mD^k(S\phi)\p_k(S\psi)- 
D^mD^k(S\psi)\p_k(S\phi) +  \\
D^k(S\phi)D^m\p_k(S\psi) - D^k(S\psi)D^m\p_k(S\phi)\big) = \\ 
D^k(S\phi)[\p_k,D^m](S\psi) - D^k(S\psi)[\p_k,D^m](S\phi) = \\
 g^{\bar l k}\p_k g^{\bar n m}\bar \p_l(S\phi)\bar \p_n(S\psi) - g^{\bar l k}\p_k g^{\bar n m}\bar \p_l(S\psi)\bar \p_n(S\phi) = \\
 g^{\bar l k}\p_k g^{\bar n m}\bar \p_l(S\phi)\bar \p_n(S\psi) - g^{\bar n k}\p_k g^{\bar l m}\bar \p_l(S\phi)\bar \p_n(S\psi) = 0, 
\end{align*}
which concludes the check that $S$ is a Poisson morphism. The proof that $T$ is an anti-Poisson morphism is similar. It remains to show that
$\{S\phi,T\psi\}_\T =0.$ It follows from (\ref{E:poisstm}), (\ref{E:st}), and (\ref{E:pspt}) that
\begin{align*}
   \{S\phi,T\psi\}_\T = \p^k S\phi\, \p_k T\psi - \bar \p^l T\psi\, \bar \p_l S\phi =  D^k S\phi \,\p_k T\psi - \\
\bar D^l T\psi \,\bar \p_l S\phi = 
g^{\bar l k} \bar \p_l S\phi\, \p_k T\psi - g^{\bar l k} \p_k T\psi\, \bar \p_l S\phi \p_k = 0.
\end{align*}
\end{proof}
According to Theorem \ref{T:main}  there exists a canonical formal symplectic groupoid ${\bf G}_U$ on the formal neighborhood $(T^*U, Z\cap U)$ such that the mappings $S,T$ are the source and target maps for ${\bf G}_U$ respectively.
The mapping $\tau: (z,\bar z, \zeta, \bar \zeta) \mapsto (z,\bar z, -\zeta, -\bar \zeta)$ is a global anti-Poisson involutive automorphism of $\T$. It induces an anti-Poisson involutive automorphism of the Poisson algebra $C^\infty(\T,Z)$. Set $\tilde S = \tau^* T$ and $\tilde T = \tau^* S$. Thus for $\phi,\psi \in C^\infty(U)$
\begin{equation}\label{E:tildest}
 (\tilde S\phi)(z, \bar z, \bar \zeta) = e^{-\bar \zeta_l \bar D^l}\phi, \ (\tilde T\psi)(z,\bar z, \zeta) = e^{-\zeta_k D^k}\psi.
\end{equation}
It follows from Proposition \ref{P:pantip} that the mappings
\[
\tilde S,\tilde T:  (C^\infty(U), \{\cdot,\cdot\}_M) \to (C^\infty(T^*U,Z\cap U),\{\cdot,\cdot\}_\T)
\]
are  a Poisson and an anti-Poisson morphisms, respectively. Moreover, for any $\phi,\psi \in C^\infty(U)$ the elements $\tilde S\phi,\tilde T\psi \in C^\infty(T^*U,Z\cap U)$ Poisson commute. Now, there is a canonical formal symplectic groupoid ${\bf \tilde G}_U$ on $(T^*U, Z\cap U)$ (the dual of ${\bf G}_U$) such that the mappings $\tilde S,\tilde T$ are the source and target maps of ${\bf \tilde G}_U$, respectively. According to formula (\ref{E:sqsprime}) there is a unique formal symplectic automorphism $Q$ of $C^\infty(T^*U,Z\cap U)$ such that 
\begin{equation}\label{E:atilde}
S = Q\tilde S \mbox{ and } T = Q\tilde T. 
\end{equation}
Let $a,\tilde a$ be arbitrary holomorphic functions and $b,\tilde b$ arbitrary antiholomorphic functions on $U$. It follows from formulas (\ref{E:st}) and (\ref{E:tildest}) that 
\begin{equation}\label{E:satb}
    Sa = a,\ Tb = b, \tilde Sb = b, \mbox{ and } \tilde Ta = a,
\end{equation}
whence we see that ${\bf G}_U$ is a formal symplectic groupoid with separation of variables over $M$ and that the dual formal groupoid ${\bf \tilde G}_U$ is a formal symplectic groupoid with separation of variables with respect to the opposite complex structure on $M$. Proposition \ref{P:pantip}, formulas (\ref{E:atilde}) and  (\ref{E:satb}) imply that
\begin{align}\label{E:qaqb}
   \{Qa,\tilde a\}_\T = \{Q\tilde Ta,\tilde a\}_\T = \{Ta,S\tilde a\}_\T=0 
   \mbox{ and}    \nonumber\\
   \{Qb, \tilde b\}_\T = \{Q\tilde Sb, \tilde b\}_\T = \{Sb, T\tilde b\}_\T=0.
\end{align}
We would like to draw the reader's attention to the analogy between formulas (\ref{E:adx}) and (\ref{E:qaqb}). There exists a unique element $F \in \J^2$ such that $Q = \exp H_F$. Represent it as 
\begin{equation}\label{E:fseries}
F = F_2 + F_3 + \ldots, 
\end{equation}
where $F_q$ is the homogeneous component of $F$ of degree $q$ with respect to the variables $\zeta_k,\bar \zeta_l$. Extracting the homogeneous components of degree $n-2$ of the left-hand sides of (\ref{E:qaqb}) and equating them to zero we obtain the following formulas where we drop the subscript $\T$ in all the Poisson brackets:
\begin{align}\label{E:faafbb}
   \left\{\left\{F_n,a\right\},\tilde a\right\} = -  \sum_{k = 2}^{n - 1} \frac{1}{k!} \sum_{i_1 + \ldots + i_k - k = n-1} \left\{\left\{F_{i_1},\ldots, \left\{F_{i_k}, a\right\} \ldots \right\}, \tilde a\right\},\nonumber \\
\left\{\left\{F_n,b\right\},\tilde b\right\} = -  \sum_{k = 2}^{n - 1} \frac{1}{k!} \sum_{i_1 + \ldots + i_k - k = n-1} \left\{\left\{F_{i_1},\ldots, \left\{F_{i_k}, b\right\} \ldots \right\}, \tilde b\right\}.
\end{align}
The right-hand sides of (\ref{E:faafbb}) depend only on $F_q$ for $q < n$ and are assumed to be equal to zero for $n=2$.
\begin{lemma}\label{L:phiphi}
  Let $\Phi_q = \Phi(z,\bar z, \zeta,\bar \zeta)$ be a homogeneous function of degree $q$ in the variables $\zeta,\bar\zeta$ on $T^*U$ such that $\{\{\Phi_q, z^i\}_\T, z^k\}_\T = 0$ and $\{\{\Phi_q, \bar z^j\}_\T, \bar z^l\}_\T = 0$ for any $i,j,k,l$. Then $\Phi_2 = \phi^{\bar l k}(z,\bar z) \zeta_k\bar \zeta_l$ for some function $\phi^{\bar l k}$ on $U$ and $\Phi_q = 0$ for $q \geq 3$.
\end{lemma}
\begin{proof}
   Using formula (\ref{E:poisstm}) we get that  $\{\{\Phi_q, z^i\}_\T, z^k\}_\T = \p^i\p^k \Phi_q = 0$ and $\{\{\Phi_q, \bar z^j\}_\T, \bar z^l\}_\T = \bar \p^j \bar \p^l \Phi_q = 0$, whence the Lemma follows. 
\end{proof}
Lemma \ref{L:phiphi} applied to formulas (\ref{E:faafbb}) implies that function (\ref{E:fseries}) is uniquely determined by the term $F_2$ which is of the form $F_2 = \phi^{\bar l k}(z,\bar z) \zeta_k\bar \zeta_l$. We can find $F_2$ explicitly using formulas (\ref{E:st}), (\ref{E:tildest}), and (\ref{E:atilde}). For an arbitrary $f = f(z,\bar z) \in C^\infty(U)$ calculate the both sides of the formula $Sf = Q(\tilde S f)$ modulo $\J^2$:
\begin{equation}\label{E:modj}
    (1 + \zeta_k D^k)f = (1 + H_{F_2})(1 - \bar \zeta_l \bar D^l) f \pmod {\J^2}.
\end{equation}
It follows from formulas (\ref{E:poisstm}) and (\ref{E:modj}) that $\p^k F_2 = g^{\bar l k}\bar \zeta_l$, whence we obtain that $\phi^{\bar l k} = g^{\bar l k}$ and therefore
\[
            F_2 = g^{\bar l k}\zeta_k \bar \zeta_l.      
\]
The remaining terms of series (\ref{E:fseries}) can be found recursively from (\ref{E:faafbb}) in local coordinates. Formula (\ref{E:even}) implies that $F_k = 0$ for the odd values of $k$.
We conclude that the function $F$ and the automorphism $Q = \exp H_F$ are uniquely determined by the K\"ahler-Poisson tensor $g^{\bar l k}$. Since condition (\ref{E:qaqb}) on $Q$ is coordinate independent, both $F$ and $Q$ are globally defined on $(\T,Z)$. It follows from formulas (\ref{E:atilde}) and (\ref{E:satb}) that for $f(z,\bar z) = a(z)b(\bar z)$
\[
    Sf = S(ab) = Sa \cdot  Sb = a \cdot Qb
\]
is completely determined by $Q$ which means that the source mapping $S$ is uniquely defined and global on $M$. The following theorem is a consequence of Theorem \ref{T:main}.
\begin{theorem}
For any K\"ahler-Poisson manifold $M$ there exists a unique formal symplectic groupoid with separation of variables on $(\T,Z)$ over $M$. Its source and target mappings are given locally by formulas (\ref{E:st}).
\end{theorem}
Now let $*$ be a star product with separation of variables on a K\"ahler-Poisson manifold $M$. Theorem \ref{T:natural} states that it is natural. The formal symplectic groupoid of the star product $*$ is the unique formal symplectic groupoid with separation of variables on $(\T,Z)$ over $M$. According to Proposition \ref{P:bertrans} the formal Berezin transform $B$ of the star product $*$ is of the form $B = \exp \frac{1}{\nu}X$, where $X$ is a natural formal differential operator on $M$. Using formula (\ref{E:commut}) we can derive from (\ref{E:adx}) and (\ref{E:qaqb}) that 
\[
\sigma(X) = F, 
\]
where $F = F_2 + F_4 + \ldots$ is determined by the condition that $F_2 = \s_2 (\Delta) = g^{\bar l k}\zeta_k \bar \zeta_l$ and equations (\ref{E:faafbb}).

\section{Appendix}

In this section we give a proof of Theorem \ref{T:ext}. To this end we need some preparations.

Let $K = (i_1,\ldots, i_n)$ be a multi-index. Denote by $K' = (i_2,i_1,\ldots,i_n)$ the multi-index obtained from $K$ by permuting $i_1$ and $i_2$, and by $\tilde K = (j_1,\ldots, j_n)$ the multi-index such that $j_1 = i_1$ and $j_2 \leq \ldots \leq j_n$ is the ordering permutation of $i_2, \ldots, i_n$. If $u^K = u^{i_1\ldots i_n}$ is a tensor symmetric in $i_2,\ldots,i_n$ then the tensor 
\begin{equation}\label{E:vuu}
v^K = u^K - u^{K'}
\end{equation}
is skew symmetric in $i_1,i_2$, symmetric in $i_3,\ldots,i_n$ and its cyclic sum over $i_1,i_2,i_3$ is zero.
\begin{lemma}\label{L:tens}
Suppose that $v^K = v^{i_1\ldots i_n}$ is a tensor skew symmetric in $i_1,i_2$, symmetric in $i_3,\ldots,i_n$ and its cyclic sum over $i_1,i_2,i_3$ is zero. There exists a unique tensor $u^K$ symmetric in $i_2,\ldots,i_n$ that satisfies (\ref{E:vuu}) and such that $u^K = 0$ if $i_1 \leq \ldots \leq i_n$.
\end{lemma}
\begin{proof}
To define $u^K$, consider $\tilde K = (j_1, \ldots, j_n)$.  Set $u^K = 0$ if $j_1 \leq j_2$ and $u^K = v^{\tilde K}$ if $j_1 \geq j_2$ (these conditions agree if $j_1 = j_2$). Thus $u^K = u^{\tilde K}$ which implies that $u^K$ is symmetric in $i_2,\ldots,i_n$. In order to show that $u^K$ is well defined we need to check condition (\ref{E:vuu}). For the multi-index $K$ in (\ref{E:vuu}) we can assume without loss of generality that $i_2 \leq i_1$ and that $i_3 = \min\{i_3,\ldots,i_n\}$. If $i_2 \leq i_3$ then $u^K = v^K$ and $u^{K'} = 0$, so (\ref{E:vuu}) holds. If $i_3 < i_2$ then $u^K = v^{i_1i_3i_2 \ldots},\ u^{K'} = v^{i_2i_3i_1\ldots}$ where the order of the remaining indices does not matter. Now (\ref{E:vuu}) holds since the cyclic sum of the tensor $v^K$ over $i_1,i_2,i_3$ is zero.
\end{proof}

For a coherent family  $\{C_n\}$ and any $f_i,\phi \in C^\infty(M)$ one can prove the following formula using Property B.
\begin{eqnarray}\label{E:phif}
  C_n(\phi, f_2, \ldots, f_n) = C_n(f_2, \ldots, f_n,\phi) + \nonumber\\
 \sum_{i = 2}^n C_{n-1}(f_2, \ldots, \{\phi, f_i\}, \ldots f_n).
\end{eqnarray}

Let $(U, \{x^i\})$ be an arbitrary coordinate chart on $M$. We will construct an operator $C_n$ locally on $U$ using induction on $n$. Assume that one can extend by one element any $k$-element coherent family for all $k < n$.
Consider an  $n$-element coherent family  $\{C_k\},  0 \leq k \leq n -1$. Then for each index $i$ and $k < n-1$ the operators 
\[
D^i_k(f_1, \ldots, f_k) = C_{k+1}(f_1,\ldots, f_k,x^i)
\]
form a coherent family. By induction this family can be extended by an operator $D^i_{n-1}$ so that 
\begin{eqnarray}\label{E:ddc}
   D^i_{n-1}(f_2, \ldots,f_k,f_{k+1}, \ldots f_n) - D^i_{n-1}(f_2, \ldots,f_{k+1},f_k, \ldots f_n)  = \nonumber \\
C_{n-1}(f_2, \ldots, \{f_k,f_{k+1}\}, \ldots, f_n,x^i),
\end{eqnarray}
Introduce the following auxilliary operator 
\begin{align}\label{E:aux}
D_n(f_1,\ldots,f_n)  =  \Big(D^i_{n-1}(f_2, \ldots,f_n) +  \nonumber\\ 
\sum_{j = 2}^n  C_{n-1}(f_2, \ldots, \{x^i, f_j\}, \ldots f_n) \Big)\frac{\p f_1}{\p x^i}.
\end{align}
The operator $D_n$ annihilates constants and is of order one in the first argument. We will show that for any $k \geq 2$
\begin{eqnarray}\label{E:cohd}
   D_n(f_1, \ldots,f_k,f_{k+1}, \ldots f_n) - D_n(f_1, \ldots,f_{k+1},f_k, \ldots f_n)  = \nonumber \\
C_{n-1}(f_1, \ldots, \{f_k,f_{k+1}\}, \ldots, f_n). 
\end{eqnarray}
Using that a derivation $A(f)$ on $U$ can be written as $A(x^i)\frac{\p f}{\p x^i}$, Property A, and formula (\ref{E:ddc}) we can show that equation (\ref{E:cohd}) is a consequence of the following one:
\begin{eqnarray}\label{E:cohdi}
 C_{n-1}(f_2, \ldots, \{f_k,f_{k+1}\}, \ldots, f_n,x^i) + \nonumber \\
\sum_{j = 2}^{k-1}
C_{n-2}(f_2, \ldots, \{x^i,f_j\},\ldots, \{f_k,f_{k+1}\}, \ldots, f_n) + \nonumber\\
\Big(C_{n-2}(f_2, \ldots, \{\{x^i, f_k\},f_{k+1}\},\ldots, f_n) + \\
C_{n-2}(f_2, \ldots, \{f_k, \{x^i, f_{k+1}\}\},\ldots, f_n)\Big) +\nonumber\\
\sum_{j = k+2}^n
C_{n-2}(f_2, \ldots, \{f_k,f_{k+1}\},\ldots,\{x^i,f_j\}, \ldots, f_n) =\nonumber\\
 C_{n-1}(x^i,f_2, \ldots, \{f_k,f_{k+1}\}, \ldots, f_n). \nonumber
\end{eqnarray}
Using the Jacobi identity, replace the sum in the parentheses in (\ref{E:cohdi}) with $C_{n-2}(f_2, \ldots, \{x^i, \{f_k, f_{k+1}\}\},\ldots, f_n)$.
The resulting identity follows from formula (\ref{E:phif}). 

We will construct the operator  $C_n$ on the coordinate chart $U$ in the form $C_n = D_n + E_n$, where $E_n(f_1,\ldots, f_n)$ is a multiderivation symmetric in $f_2, \ldots, f_n$. The operator $E_n$ must be chosen so that $C_n$ would satisfy Property B for $k =1$ (all other conditions on $C_n$ are already satisfied). 
This condition can be written in the form
\begin{equation}\label{E:iv}
     V_n(f_1,f_2,\ldots, f_n) = E_n(f_2,f_1,\ldots,f_n) - E_n(f_1,f_2,\ldots,f_n),
\end{equation}
where the operator $V_n$ is given by the formula
\begin{align}\label{E:vn}
   V_n(f_1,f_2,\ldots, f_n) = D_n(f_1,f_2,\ldots, f_n) - 
D_n(f_2,f_1,\ldots, f_n) \nonumber\\ 
- C_{n-1}(\{f_1,f_2\},\ldots, f_n).  
\end{align}
According to Lemma \ref{L:tens}, an operator $E_n$ with the required properties exists if $V_n(f_1,f_2,\ldots, f_n)$ is a multiderivation skew symmetric in $f_1,f_2$, symmetric in $f_3,\ldots,f_n$, and such that the cyclic sum of $V_n$ over $f_1,f_2,f_3$ is zero. We will show that the operator $V_n$ enjoys all these properties. Check that the operator $V_n$ is a derivation in the second argument. Substituting formula (\ref{E:aux}) in (\ref{E:vn}) and taking into account Property A we see that it remains to check that the operator
\begin{equation}\label{E:rem}
       C_{n-1}(\{x^i,f_2\},f_3,\ldots,f_n)\frac{\p f_1}{\p x^i} -  C_{n-1}(\{f_1,f_2\},\ldots, f_n) 
\end{equation}
is a derivation in $f_2$. Formula (\ref{E:rem}) can be rewritten as follows,
\begin{equation}\label{E:exp}
   C_{n-1}(x^j,f_3,\ldots,f_n)\left(\frac{\p f_1}{\p x^i} \frac{\p}{\p x^j}\{x^i,f_2\} -  \frac{\p}{\p x^j}\{f_1,f_2\}\right).
\end{equation}
In local coordinates $\{f,g\} = \eta^{kl}\frac{\p f}{\p x^k}\frac{\p g}{\p x^l}$, where $\eta^{kl}$ is a Poisson tensor. The second factor in (\ref{E:exp}) equals
\[
  \frac{\p f_1}{\p x^k} \frac{\p}{\p x^j}\left(\eta^{kl} \frac{\p f_2}{\p x^l}\right) -    \frac{\p}{\p x^j} \left(\eta^{kl} \frac{\p f_1}{\p x^k}\frac{\p f_2}{\p x^l}\right)  =
- \eta^{kl} \frac{\p^2 f_1}{\p x^j \p x^k}\frac{\p f_2}{\p x^l}.
\]
Thus $V_n(f_1,f_2,\ldots,f_n)$ is a derivation in $f_2$. Since it is obviously skew symmetric in $f_1,f_2$, it is also a derivation in $f_1$.

We will prove that $V_n(f_1,f_2,\ldots,f_n)$ is symmetric in $f_3, \ldots,f_n$ using formula (\ref{E:cohd}). For $k \geq 3$
\begin{align*}
   V_n(f_1,\ldots,f_k,f_{k+1}, \ldots, f_n) -  V_n(f_1,\ldots,f_{k+1},f_k, \ldots, f_n) = \\
C_{n-1}(f_1,f_2,\ldots,\{f_k,f_{k+1}\},\ldots,f_n) - \\ C_{n-1}(f_2,f_1,\ldots,\{f_k,f_{k+1}\},\ldots,f_n) -\\ C_{n-2}(\{f_1,f_2\},\ldots,\{f_k,f_{k+1}\},\ldots,f_n) =0.
\end{align*}
It remains to show that $V_n(f_1,f_2,\ldots,f_n)$ is a derivation in $f_3$ and that its cyclic sum over $f_1,f_2,f_3$ is zero.  We have, using formula (\ref{E:cohd}), that
\begin{align*}
  V_n(f_1,f_2,f_3, \ldots, f_n) = D_n(f_1,f_2,f_3,\ldots,f_n) - \\D_n(f_2,f_1,f_3,\ldots, f_n) - C_{n-1}(\{f_1,f_2\},f_3,\ldots, f_n) =\\
D_n(f_1,f_3,f_2,\ldots, f_n) + C_{n-1}(f_1,\{f_2,f_3\}, \ldots, f_n) - \\ D_n(f_2,f_1,f_3, \ldots, f_n) - C_{n-1}(\{f_1,f_2\},f_3, \ldots, f_n ) = \\
D_n(f_1,f_3,f_2, \ldots, f_n) - D_n(f_2,f_1,f_3, \ldots, f_n) + \\C_{n-2}(\{f_1,\{f_2,f_3\}\}, \ldots, f_n).
\end{align*}
We see that the cyclic sum of $V_n$ over $f_1,f_2,f_3$ is zero due to the Jacobi identity. Therefore,
\begin{align}\label{E:vvv}
   V_n(f_1,f_2,f_3, \ldots, f_n) = - V_n(f_2,f_3,f_1, \ldots, f_n) - \nonumber\\V_n(f_3,f_1,f_2, \ldots, f_n).
\end{align}
We have already proved that $V_n$ is a derivation in the first two arguments, whence the right hand side and therefore the left hand side of (\ref{E:vvv}) are derivations in $f_3$. Since $V_n(f_1,\ldots,f_n)$ was shown to be symmetric in $f_3,\ldots,f_n$, this implies that $V_n$ is a multiderivation. This concludes the proof of all the properties of the operator $V_n$ and provides a local construction of the operator $C_n$. 
Finally, we use partition of unity to construct $C_n$ globally on $M$.

\end{document}